\documentclass[12pt]{amsart}

\numberwithin{equation}{section}
\def\a{\alpha}

\def\b{\beta}

\def\o{\omega}

\def\spt{\text{spt}\,}

\def\span{\text{span}\,}
\def\th{\theta}

\newcommand\trace{\operatorname{trace}}
\newcommand\w{\wedge}

\newcommand\R{\mathbf R}
\newcommand\rp{\mathbf {RP}}
\newcommand\Z{\mathbf Z}

\def\1{\{1\}}
\def\2{\{2\}}
\def\3{\{3\}}
\def\4{\{4\}}

\newtheorem*{introtheorem1}{Theorem \ref{MT3}}
\newtheorem*{introtheorem2}{Theorem \ref{5T2}}
\newtheorem{theorem}{Theorem}
\newtheorem{lemma}[theorem]{Lemma}
\newtheorem{corollary}[theorem]{Corollary}
\newtheorem{proposition}[theorem]{Proposition}

\theoremstyle{definition}
\newtheorem{definition}[theorem]{Definition}

\theoremstyle{remark}
\newtheorem{remark}[theorem]{Remark}

\begin{document}

\title[Dupin hypersurfaces]{Dupin hypersurfaces with four principal
curvatures, II} 

\author[Cecil]{Thomas E. Cecil}
\thanks{The first author was partially supported by NSF Grant
No. DMS-0405529} 
\address{Department of Mathematics and Computer Science \\ College 
of the Holy Cross \\ Worcester, Massachusetts 01610-2395} 
\email{cecil@mathcs.holycross.edu}
\author[Chi]{Quo-Shin Chi}
\address{Department of Mathematics \\ Campus Box 1146 \\ Washington
University \\ St. Louis, Missouri 63130} 
\email{chi@math.wustl.edu}
\author[Jensen]{Gary R. Jensen}
\address{Department of Mathematics \\ Campus Box 1146 \\ Washington
University \\ St. Louis, Missouri 63130} 
\email{gary@math.wustl.edu}

\date{\today}


\begin{abstract} If $M$ is an isoparametric hypersurface in a sphere
$S^n$ with four distrinct principal curvatures, then the principal
curvatures $\kappa_1,\dots,\kappa_4$ can be ordered so that their
multiplicities satisfy $m_1=m_2$ and $m_3=m_4$, and the cross-ratio
$r$ of the principal curvatures (the Lie curvature) equals $-1$.  In
this paper, we prove that if $M$ is an irreducible connected proper
Dupin hypersurface in $\R^n$ ( or $S^n$) with four distinct principal
curvatures with multiplicities $m_1=m_2 \geq 1$ and $m_3=m_4=1$, and
constant Lie curvature $r=-1$, then $M$ is equivalent by Lie sphere
transformation to an isoparametric hypersurface in a sphere.  This
result remains true if the assumption of irreducibility is replaced
by compactness and $r$ is merely assumed to be constant.
\end{abstract}

\maketitle


\section{Introduction}
Let $M$ be an immersed hypersurface in Euclidean space ${\bf R}^n$
or the unit sphere $S^n \subset {\bf R}^{n+1}$. A {\em curvature surface} 
of $M$ is a smooth connected submanifold $S$
such that for each point $x \in S$, the tangent space
$T_xS$ is equal to a principal space of the shape operator
$A$ of $M$ at $x$.  This generalizes the classical notion of a 
line of curvature on a surface in ${\bf R}^3$.  The hypersurface $M$
is said to be \textit{Dupin} if it satisfies the condition

\begin{enumerate}
\item[(a)] along each curvature surface, the corresponding principal
curvature is constant.
\end{enumerate}
The hypersurface $M$ is called {\em proper Dupin} if, in addition 
to condition (a), it also satisfies the condition

\begin{enumerate}
\item[(b)] the number $g$ of distinct principal curvatures is constant
on $M$.
\end{enumerate}
Pinkall \cite{Pinkall} proved that
both of these conditions are invariant under the group of
Lie sphere transformations of $S^n$, which contains the
group of M\"{o}bius (conformal) transformations of $S^n$ as
a subgroup.  Thus, by stereographic projection, the theory of
Dupin hypersurfaces in ${\bf R}^n$ or $S^n$ is essentially the same.

Thorbergsson \cite{Th1} showed that the number $g$ of distinct 
principal curvatures of a compact proper Dupin hypersurface $M$
immersed in $S^n$ must be $1,2,3,4$ or 6,
the same as M\"{u}nzner's \cite{Mun1, Mun2} restriction on
the number of distinct  
principal curvatures of an isoparametric (constant principal curvatures)
hypersurface in $S^n$.
In the cases $g = 1,2,3$, compact proper Dupin hypersurfaces
in $S^n$ have been completely classified.  For the totally umbilic case
$g=1$, $M$ must be a great or small sphere.  For $g=2$, Cecil and Ryan
\cite{CecilRyan2} 
proved that $M$ must be M\"{o}bius equivalent to a standard product
of spheres (which is isoparametric)
\begin{displaymath}
S^k(r) \times S^{n-k-1}(s) \subset S^n, \quad r^2+s^2=1.
\end{displaymath}
In the case $g=3$, Miyaoka \cite{Mi1} proved that $M$ must be Lie
equivalent 
to an isoparametric hypersurface in $S^n$, which by the work of
Cartan \cite{Car2} must be a tube of
constant radius over a standard Veronese embedding of a projective
plane $FP^2$ into $S^{3m+1}$, where $F$ is the division algebra
{\bf R}, {\bf C}, {\bf H} (quaternions),
{\bf O} (Cayley numbers) for $m=1,2,4,8,$ respectively.

The cases of compact proper Dupin hypersurfaces with $g = 4$ or 6
principal 
curvatures have not yet been classified, although Stolz \cite{Stolz}
in the case $g=4$ 
and Grove and Halperin \cite{GH} in the case $g=6$ have shown that the
multiplicities 
of the principal curvatures of a compact proper Dupin hypersurface
must be  
the same as for an isoparametric hypersurface.  
In particular, in the case 
$g=4$, the multiplicities must come in pairs, and the principal
curvatures can 
be ordered in such a way that $m_1 = m_2$ and $m_3 = m_4$.  In the
case $g=6$, 
all the principal curvatures must have the same multiplicity $m = 1$
or 2. 

Miyaoka \cite{Mi2} introduced an important set of 
Lie invariants,  the {\it Lie curvatures} of a Dupin hypersurface $M$,
which are the  
cross-ratios of the principal curvatures taken four at a time.
Obviously,  
for an isoparametric hypersurface, the Lie curvatures are constant,
and a 
necessary condition for a Dupin hypersurface to be Lie equivalent to an
isoparametric hypersurface is that it have constant Lie curvatures.
At one time it was thought that perhaps every compact proper Dupin is
Lie equivalent 
to an isoparametric hypersurface.  However,
by two separate constructions, Pinkall and Thorbergsson \cite{PT}
$(g=4)$ and  
Miyaoka and Ozawa \cite{MiOz} ($g = 4$ or 6) produced compact proper
Dupin  
hypersurfaces which do not have constant Lie curvatures
and therefore cannot be Lie equivalent to an isoparametric
hypersurface. 

Miyaoka \cite{Mi2,Mi3} showed that a compact proper Dupin hypersurface 
immersed in $S^n$ with $g = 4$ or $6$ principal curvatures 
is Lie equivalent to an isoparametric hypersurface if it has constant
Lie curvatures and it satisfies certain 
additional global conditions regarding the intersections of leaves of
its various principal foliations.  The goal of current research is to
prove  
that the condition of constant Lie curvatures already suffices for the
conclusion 
without assuming these additional conditions, and we have succeeded 
in doing this in the case $g=4$ when one pair of the multiplicities is
equal to one, as described below.

In contrast to the situation for compact proper Dupin hypersurfaces,
there is a local method, due to Pinkall \cite{Pinkall}, for producing
a Dupin hypersurface with any given number $g$ of
principal curvatures with any prescribed multiplicities
$m_1,\ldots,m_g$.  His method uses 
the basic constructions of building tubes,
cylinders, cones and surfaces of revolution over a
Dupin hypersurface
$W^{n-1}$ in ${\bf R}^n$ with $g$ principal curvatures to
get a Dupin hypersurface $M^{n-1+k}$ in
${\bf R}^{n+k}$ with $g+1$ principal curvatures.
These constructions introduce a new principal curvature 
of multiplicity $k$ which is easily seen to be constant along its 
curvature surfaces.  The other principal curvatures are determined by the 
principal curvatures of $W^{n-1}$, and the Dupin property is preserved
for these principal curvatures. These constructions
are local in nature and
only yield a compact proper Dupin hypersurface if the original
manifold $W^{n-1}$ is itself a sphere \cite[Theorem 46]{Cec5}.  Otherwise, the
number of 
distinct principal curvatures is not constant on a compact manifold
$M^{n-1+k}$ obtained in this way, so it is not proper Dupin.

A Dupin hypersurface which is locally equivalent 
by a Lie sphere transformation to a
hypersurface 
$M^n$ obtained by one of these constructions is said to be 
{\it reducible}.  Otherwise, the Dupin hypersurface is
called {\it irreducible}.  A Dupin hypersurface is called
{\it locally irreducible} if it does not contain any reducible
open subset.  Clearly, local irreducibility implies irreducibility.

In~\cite{AlgDupin}, we prove that any $C^\infty$ proper
Dupin hypersurface  
must be analytic.  Using
this result, we 
prove that if a connected proper Dupin hypersurface $M$ has a
reducible open subset, 
then $M$ itself is reducible.  That is, irreducibility implies
local irreducibility. Analyticity allows us to work locally to obtain
global results.

The primary work in this paper is local in nature and is
accomplished in the setting of Lie sphere geometry.  
We concentrate 
on the case $g=4$ with multiplicities $m_1 = m_2$, $m_3 = m_4$ and 
Lie curvature $r = -1$.  In Section~\ref{Lie} we review the concepts
of Lie sphere geometry  as well as
the basic set-up of the method of moving frames developed in our
previous  
paper \cite{CJ2}.  

In Section~\ref{Isopara} we specialize the Lie frame for the case
under study.  Theorem~\ref{Th11} establishes sufficient conditions
for when a proper Dupin hypersurface with four principal curvatures
having multiplicities $m_1 = m_2$,
$m_3 = m_4$ and  
Lie curvature $r = -1$ is Lie equivalent to an isoparametric
hypersurface. 

In Section~\ref{Necessary}, we relate the Lie
sphere definition of reducibility to the Pinkall constructions in
Euclidean space.  
Theorems~\ref{TS2} and~\ref{PU1} establish sufficient  
conditions for reducibility in terms of quantities that arise
naturally in the setting of moving Lie frames.

In Section~\ref{Sec8}, we prove one of our main
results:

\begin{introtheorem1} Suppose the connected proper Dupin
hypersurface $\lambda: M^{n-1} \to \Lambda^{2n-1}$ has four distinct
curvature spheres with multiplicities
$m_1=m_2 \geq 1$, $m_3=m_4=1$, and Lie curvature $r = -1$.  If
$\lambda$ is irreducible, then it is Lie equivalent to an
isoparametric hypersurface.
\end{introtheorem1}
 
In~\cite[p.3]{CJ2}, it was conjectured that if $M$ is an irreducible
proper Dupin hypersurface in $S^n$ with four principal curvatures
having respective multiplicities $m_1$, $m_2$, $m_3$, $m_4$, and $M$
has constant Lie curvature, then the principal curvatures can be
ordered so that $m_1=m_2$, $m_3=m_4$, and $M$ is Lie equivalent to an
isoparametric hypersurface in $S^n$.  We still believe this conjecture
to be true, although we have not yet been able to verify it in more
generality than Theorem~\ref{MT3}.

In Section~\ref{Compact}, we prove in Theorem~\ref{5T1} that a compact
proper Dupin hypersurface with $g >2$ principal curvatures is
irreducible.  As a consequence of this, Theorem~\ref{MT3}, and a
result of Miyaoka \cite{Mi2} that a compact proper Dupin hypersurface
with $g=4$ and constant Lie curvature $r$ must have $r=-1$, we obtain
our second main result:

\begin{introtheorem2} Let $M$ be a compact connected proper Dupin
hypersurface immersed in $R^n$ with four distinct principal curvatures
having multiplicities $m_1=m_2 \geq 1$, $m_3 = m_4 = 1$, and 
constant Lie curvature.  Then $M$
is Lie equivalent to an isoparametric hypersurface.
\end{introtheorem2}


\section{Dupin hypersurfaces in Lie sphere geometry}\label{Lie}

In this section, we briefly recall how Dupin hypersurfaces
can be studied in the context of Lie sphere geometry.  In
particular, we will summarize the basic set-up and
main definitions of \cite{CJ2} that will be needed in the remainder
of the paper.  We will not, however, reproduce all the formulas
from that paper, so the reader will need
to consult that paper at times.  Throughout this paper, equation
references of the sort GD(3.36) will be to equation (3.36) of
\cite{CJ2}.  We will use
the Einstein summation convention in this section.
  
Let $\R_2^{n+3}$ be a real vector space of
dimension $n+3$ endowed with the metric of signature $(n+1,2)$,
\begin{equation}\label{Tom1}
\langle x,y \rangle = -x^0 y^0 + x^1 y^1 + \cdots +x^{n+1} y^{n+1} -
x^{n+2} y^{n+2}.
\end{equation}
Let $e_0, \ldots ,e_{n+2}$ denote the standard orthonormal
basis with respect to this metric, with $e_0$ and $e_{n+2}$
timelike.  Let $P^{n+2}$ be the real projective space of lines through
the origin in $\R_2^{n+3}$,  and let $Q^{n+1}$ be the quadric
hypersurface determined by the equation $\langle x,x \rangle = 0$.  This
hypersurface is called the {\it Lie quadric}.
The sphere $S^n$ can be identified with the unit sphere in the Euclidean space
$\R^{n+1}$ spanned by the vectors $e_1, \ldots , e_{n+1}$.

The points in $Q^{n+1}$ are in bijective correspondence
with the set of all oriented hyperspheres and point
spheres in $S^n$.  
The Lie quadric contains projective lines but no linear
subspaces of $P^{n+2}$ of higher dimension. 
Let $\Lambda^{2n-1}$ denote the set of all projective lines in
$Q^{n+1}$.  It is an analytic manifold of dimension $2n-1$.
The line $[x, y]$
determined by two points $[x]$ and $[y]$ of $Q^{n+1}$ lies on
$Q^{n+1}$ if and only if $\langle x,y \rangle = 0.$  This happens
precisely when the hyperspheres in $S^n$ corresponding to the
points $[x]$ and $[y]$ are in oriented contact.

A \textit{Lie sphere transformation} is a projective transformation
of $P^{n+2}$ which maps $Q^{n+1}$ to itself.  A Lie sphere
transformation 
preserves oriented contact of hyperspheres in $S^n$, since it takes
lines on $Q^{n+1}$ to lines on $Q^{n+1}$.  The group $G$ of Lie sphere
transformations is isomorphic to $O(n+1,2)/\{\pm I\},$ where
$O(n+1,2)$ is the orthogonal group for the metric~\eqref{Tom1}.  The
group $G$ acts transitively on $\Lambda^{2n-1}$.
  
The manifold  $\Lambda ^{2n-1}$ of projective lines on $Q^{n+1}$
has a \textit{contact structure},
i.e., a globally defined 1-form $\omega$ such that
$\omega \wedge d\omega ^{n-1}$ never vanishes on
$\Lambda ^{2n-1}$.  The condition $\omega = 0$ defines a
codimension one distribution $D$ on $\Lambda ^{2n-1}$ which has
integral submanifolds of dimension $n-1$ but none of higher
dimension.  A \textit{Legendre submanifold} is one of these integral
submanifolds of maximal dimension, i.e.,
an immersion $\lambda : M^{n-1} \rightarrow
\Lambda ^{2n-1}$ such that $\lambda ^* \omega = 0$.

An immersion $f:M^{n-1} \rightarrow S^n$ with field of unit
normals $\xi :M^{n-1} \rightarrow S^n$ naturally induces
a Legendre submanifold $\lambda = [Y_0, Y_1]$, where
$[Y_0, Y_1]$ denotes 
the line in $Q$ determined by the point sphere
$Y_0 = (1,f,0)$ and the tangent great sphere  $Y_1 = (0,\xi ,1)$.
In a similar way,
an immersed submanifold $\phi :V \rightarrow S^n$ of codimension 
greater than one also induces a Legendre submanifold 
whose domain is the bundle $B^{n-1}$ of unit normal
vectors to $\phi (V)$ (see, for example, \cite[p.79]{Cec6}). 

Suppose that $\lambda = [Y_0 , Y_1]$ is a Legendre submanifold.  Let
$p \in M^{n-1}$ and let $r$ and $s$ be real numbers at least one
of which is non-zero.  The sphere in $S^n$ corresponding to the point
$[K] = [rY_0 (p) + sY_1 (p) ]$ in $Q^{n+1}$
is called a {\it curvature sphere} of $\lambda$ at $p$, if there
exists a non-zero tangent vector $X \in T_p M$ such that
$r dY_0 (X) + s dY_1 (X) \in {\rm Span} \{ Y_0 (p),
Y_1 (p) \}$. The vector $X$ is called a {\it principal vector}
corresponding to the curvature sphere $[K]$.  The principal
vectors corresponding to a given curvature sphere form a subspace
of $T_p M$, and $T_pM$ is the direct sum of these {\it principal
spaces}. 

To see the relationship between curvature spheres and principal
curvatures, suppose now that
$\lambda = [Y_0 , Y_1 ]$ with $Y_0 = (1,f,0),  Y_1 = (0,\xi ,1)$.  
At a given
$p \in M$, one can write the distinct curvature spheres
in the form
$[K_i ] = [\kappa _i Y_0 + Y_1 ], 1\leq i \leq g.$
In the case where the map $f$ is an immersion, these
$\kappa _i$ are the usual principal curvatures
of the hypersurface $f$ at $p$. 
The principal
curvatures are not invariant
under Lie sphere transformations.  However, the cross-ratio
of any four distinct principal curvatures is Lie invariant.
These cross ratios are called {\it Lie curvatures} of $\lambda$.

As in Euclidean submanifold theory, a curvature surface is
a smooth connected submanifold $S$ of $M$ such that for each point
$p \in S$, the tangent space $T_p S$ is equal to a principal
space.  A Legendre submanifold is called \textit{Dupin} if
along each curvature surface, the corresponding curvature sphere
is constant.  A Dupin submanifold is said to be \textit{proper Dupin}
if the number $g$ of distinct curvature spheres is constant on $M$.
These definitions agree with the usual definitions in the case
where the Legendre submanifold is induced from an immersed
hypersurface in $S^n$.  Pinkall \cite{Pinkall} showed that both of these
properties are invariant under the group of Lie sphere transformations.
At times, we will refer to Dupin submanifolds as ``Dupin hypersurfaces,''
because of their close relationship with Dupin hypersurfaces in $S^n$.

We now begin to recall the notation and results from \cite{CJ2} in detail.
We study Dupin hypersurfaces in Lie sphere geometry using the
method of moving frames.
Instead of using an orthonormal frame for the metric in~\eqref{Tom1}, we consider
a \textit{Lie frame}, that is, an ordered set of vectors
$Y_0,\ldots,Y_{n+2}$ in $\R_2^{n+3}$ satisfying
$\langle Y_a,Y_b \rangle = k_{ab}$, for $0 \leq a,b \leq n+2$,
with
\begin{equation}\label{GD(2.4)}
k = (k_{ab}) = \begin{pmatrix} 0&0& -J \\ 0&I_{n-1} &0\\ -J &0&0 \end{pmatrix}, \quad
\text{ where } \quad J= \begin{pmatrix} 0&1\\1&0 \end{pmatrix}
\end{equation}
The space of all Lie frames can be identified with the
orthogonal group $G = O(n+1,2).$
In this space, one introduces the Maurer-Cartan forms,
\begin{equation}
dY_a = \omega^b_a Y_b, 0 \leq a,b \leq n+2,
\end{equation}
which satisfy the Maurer-Cartan structure equations of $G$,
\begin{equation}\label{GD(2.21)}
d\o^a_b = - \o^a_c \w \o^c_b, \quad \text{ for $0 \leq a,b,c \leq
n+2$}
\end{equation}
Knowing from \cite{AlgDupin} that any proper Dupin hypersurface is
real analytic, we assume from now on that all maps are real analtyic.
A {\it Lie frame field} along a Legendre submanifold
$\lambda : M^{n-1} \rightarrow
\Lambda ^{2n-1}$ is a real analytic map $Y:U \rightarrow G$ defined
on an open subset $U$ of $M^{n-1}$
such that $\lambda(p) = [Y_0(p),Y_1(p)]$ for each $p\in U$.
Here $Y_a$ denotes the $a^{\text{th}}$ column of $Y$ and $Y \in G$
means $\langle Y_a, Y_b \rangle = k_{ab}$, for all $a,b =
0,1,\dots,n+2$. 

The notion of a curvature sphere of a Legendre submanifold 
$\lambda : M^{n-1} \to \Lambda ^{2n-1}$ can be formulated in 
terms of Lie frames as follows.  If $Y$ is any Lie frame
field along $\lambda$ defined on a neighborhood of a point $p \in M$,
then $[rY_0 + sY_1]$ is a curvature sphere of $\lambda$ at $p$ 
precisely when the following equation is satisfied at $p$,
\begin{equation}\label{GD(2.37)}
(r\o^2_0 +s \o^2_1)\w \dots \w (r\o^n_0 + s\o^n_1) = 0
\end{equation}
This condition is equivalent to saying that the tangent sphere map
\begin{equation}\label{2.37a}
[rY_0 + sY_1]: U \to Q \subset P^{n+2}
\end{equation}
is singular at $p$ in the sense that there exists a non-zero vector $X
\in T_pM$ such that
\begin{equation}\label{2.37b}
d(rY_0 + sY_1)_{(p)}(X) \in \span \{Y_0(p), Y_1(p)\}
\end{equation}
We now restrict our attention to the case where the Legendre
submanifold 
$\lambda : M^{n-1} \rightarrow \Lambda ^{2n-1}$ has $g=4$ distinct
curvature spheres of multiplicities $m_1$, $m_2$, $m_3$ and $m_4$,
respectively. We define sets
\begin{equation}\label{GD(2.44)}
\aligned
\1 &= \{2, \dots, m_1+1\} \\  
\2 &= \{m_1+2, \dots, m_1+m_2+1\} \\ 
\3 &= \{m_1+m_2+2, \dots, m_1 + m_2 + m_3 +1\} \\
\4 &= \{m_1+m_2+m_3+2,\dots, m_1+m_2+m_3+m_4+1\}.
\endaligned
\end{equation}
and adopt the index conventions
\begin{equation}\label{GD(2.20_2.45)}
\aligned
2 \leq i,j,k,l & \leq n \\
a,b,c,d &\in \1 \\
p,q,r,s &\in \2 \\
\a,\b,\gamma,\delta &\in \3 \\
\mu,m,\sigma,\tau &\in \4
\endaligned
\end{equation}
We next recall the following definition from \cite{CJ2}.

\begin{definition}\label{FirstOrderDef}
Suppose that $\lambda :M\to \Lambda$ 
is a real analytic Legendre submanifold with
$g=4$ distinct curvature spheres at each point.
A \textit{first order frame field along} $\lambda$ is an analytic
Lie frame field $Y : U \subset M\to G$ such that
\begin{equation}\label{GD(2.46)}
[Y_0], \quad [Y_1], \quad [Y_0 + Y_1], \quad [rY_0+Y_1]
\end{equation}
are the curvature spheres of
$\lambda$ at each point of $U$, and
\begin{equation}\label{GD(2.47)}
\o^a_0 =0, \quad \o^p_1 =0, \quad \o^\a_0 + \o^\a_1 =0, \quad
r\o^\mu_0+\o^\mu_1=0
\end{equation}
for all $a,p,\a,\mu$.
\end{definition}
Here $r:M \rightarrow \R$ is an analytic function never taking the
values 0 or 1.  Since we are free to put the four curvature spheres
in any order, we can assume
\begin{equation}\label{GD(2.48)}
-\infty < r <0
\end{equation}
Note that $r$ is the cross-ratio of the curvature spheres in the
appropriate order, and thus $r$ is the Lie curvature of $\lambda$.

One can show that there exists a first order Lie frame field defined
on some neighborhood of any point of any analytic Legendre
submanifold with $g=4$ distinct curvature spheres at each point.
If $Y$ is a first order frame field on
an open set $U \subset M$, then its \textit{associated coframe field}
in $U$ is the set of analytic 1-forms
\begin{equation}\label{GD(2.54)}
\theta^a = \o^a_1, \quad\theta^p = \o^p_0, \quad\theta^\a = \o^\a_0,
\quad\theta^\mu = \o^\mu_0
\end{equation}
We now assume that the Legendre submanifold $\lambda : M^{n-1}
\rightarrow \Lambda ^{2n-1}$ is connected and proper Dupin with constant Lie
curvature $r$. 

\begin{definition}\label{SecondOrder} A \textit{second-order Lie frame
field along} $\lambda$ is a  
first order frame field $Y:U \to G$ such that
\begin{equation}\label{GD(3.12)}
\o^1_0 = 0, \quad \o^0_1 = 0, \quad \o^1_1 = \o^0_0
\end{equation}
\end{definition}
In \cite{CJ2}, we show that it follows from the Dupin condition,
that for any point $p \in M$, there exists a neighborhood $U$ of $p$
on which there is defined a second-order frame field along $\lambda$.
If $Y:U \to G$ is a second order Lie frame field on an open set $U
\subset M$, then any other second order frame field on $U$ is given by
\begin{equation}\label{GD(3.14)}
\tilde Y = Y\, a(tI_2,B,0,sL)
\end{equation}
where $t,s:U \to \R$ are real analytic functions, $t$ never zero,
\begin{equation}\label{GD(3.15)}
a(tI_2,B,0,sL) = \begin{pmatrix} tI_2 & 0 & sL \\
0 & B & 0 \\
0 & 0 & t^{-1}I_2 \end{pmatrix}
\end{equation}
\begin{equation}\label{LL}
L = \begin{pmatrix} 1&0 \\ 0&-1 \end{pmatrix}
\end{equation}
and
\begin{equation}\label{BB}
B = \begin{pmatrix} B_1 &0 & 0 & 0 \\
0& B_2 &0&0 \\
0&0& B_3 &0 \\
0&0&0 & B_4 \end{pmatrix}
\end{equation}
where $B_i:U \to O(m_i)$ are real analytic maps.

With a second order frame, we have the following basic expressions
for certain Maurer-Cartan forms in terms of the associated coframe.
These equations define the analytic tensors $F^\alpha_{pa}$,
$F^\mu_{pa}$, $F^\mu_{\alpha a}$, and $F^\mu_{\alpha p}$, which
are crucial to our study.
\begin{equation}\label{GD(3.13)}
\aligned
\o^p_a &= F^\alpha_{pa}\th^\alpha + rF^\mu_{pa}\th^\mu \\
\o^\alpha_a &= F^\alpha_{pa}\th^p + (r-1)F^\mu_{\alpha a}\th^\mu \\
\o^\alpha_p &= F^\alpha_{pa}\th^a + (r-1)F^\mu_{\alpha p} \th^\mu \\
\o^\mu_a &= rF^\mu_{pa}\th^p + (r-1)F^\mu_{\alpha a}\th^\alpha \\
\o^\mu_p &= F^\mu_{pa}\th^a + \frac{r-1}r F^\mu_{\alpha p}\th^\alpha \\
\o^\mu_\alpha &= F^\mu_{\alpha a}\th^a + F^\mu_{\alpha p}\th^p
\endaligned
\end{equation}
We also have the following formulas for the Maurer-Cartan forms 
$\omega^0_i$ and $\omega^1_i$,
\begin{equation}\label{GD(3.24)}
\o^0_i = D_{ij}\th^j \text{ and } \o^1_i=E_{ij}\th^j
\end{equation}
The conditions defining a second order frame together with the
structure equations impose many conditions on the analytic functions
$D_{ij}$ 
and $E_{ij}$, which are listed in equations GD(3.25) and GD(3.26)
of~\cite{CJ2}. In summary, what emerges are
eight symmetric matrices
\begin{equation}\label{GD(3.27)}
\alignedat4
D_1 &= (D_{ab}) \quad& D_2 &= (D_{pq})\quad& D_3 &= (D_{\alpha\beta})
\quad& D_4 &= (D_{\mu\nu}) \\ 
E_1 &= (E_{ab}) & E_2 &= (E_{pq}) & E_3 &= (E_{\alpha\beta}) & E_4 &=
(E_{\mu\nu}) 
\endalignedat
\end{equation}
and six matrices of analytic functions
\begin{equation}\label{GD(3.28)}
D_{a\alpha}, \quad D_{pa}, \quad D_{p\alpha}, \quad D_{\mu a},\quad
D_{\mu\alpha}, \quad E_{p\mu}
\end{equation}
so that equations~\eqref{GD(3.24)} become
\begin{equation}\label{GD(3.29)}
\aligned
\o^0_a &= D_{ab}\th^b + D_{a\alpha} \th^\alpha - rD_{\mu a}\th^\mu \\
\o^0_p &= D_{pa}\th^a + D_{pq}\th^q + D_{p\alpha}\th^\alpha +
rE_{p\mu}\th^\mu \\ 
\o^0_\alpha &= -D_{a\alpha}\th^a + D_{\alpha\beta}\th^\beta +
rD_{\mu\alpha} \th^\mu \\ 
\o^0_\mu &= D_{\mu a}\th^a + D_{\mu\alpha}\th^\alpha +
D_{\mu\nu}\th^\nu 
\endaligned
\end{equation}
and
\begin{equation}\label{GD(3.30)}
\aligned
\o^1_a &= E_{ab}\th^b - D_{pa}\th^p + D_{a\alpha}\th^\alpha - D_{\mu
a}\th^\mu \\ 
\o^1_p &= E_{pq}\th^q + D_{p\alpha}\th^\alpha + E_{p\mu}\th^\mu \\
\o^1_\alpha &= D_{p\alpha}\th^p + E_{\alpha\beta}\th^\beta +
D_{\mu\alpha}\th^\mu \\ 
\o^1_\mu &= E_{p\mu}\th^p + D_{\mu\alpha}\th^\alpha +
E_{\mu\nu}\th^\nu 
\endaligned
\end{equation}
The tensors in~\eqref{GD(3.27)} satisfy the set of four linear
equations GD(3.42), which relate these functions to the four
multiplicities. 
The functions in~\eqref{GD(3.28)} also arise in the following
important expression for the exterior derivative of the 
form $\omega^0_0$,
\begin{equation}\label{GD(3.33)}
\aligned
d\o^0_0 &= -D_{pa}\th^a\w\th^p + D_{a\alpha}\th^a\w\th^\alpha - D_{\mu
a}\th^a\w\th^\mu \\ 
&+ D_{p\alpha} \th^p\w\th^\alpha + rE_{p\mu} \th^p\w\th^\mu + (r-1)
D_{\mu\alpha} \th^\alpha\w\th^\mu 
\endaligned
\end{equation}
In this set-up, we define the covariant derivatives
of the $F$'s, as  the analytic functions on
the right side of the equations
\begin{equation}\label{GD(3.34)}
\aligned
dF^\alpha_{pa} + F^\alpha_{pa}\o^0_0 + F^\beta_{pa}\o^\alpha_\beta -
F^\alpha_{pb}\o^b_a - F^\alpha_{qa}\o^q_p &= F^\alpha_{paj}\th^j \\ 
dF^\mu_{pa} + F^\mu_{pa}\o^0_0 + F^\nu_{pa}\o^\mu_\nu -
F^\mu_{pb}\o^b_a - F^\mu_{qa}\o^q_p &= F^\mu_{paj}\th^j \\ 
dF^\mu_{\alpha a} + F^\mu_{\alpha a}\o^0_0 + F^\nu_{\alpha
a}\o^\mu_\nu - F^\mu_{\alpha b}\o^b_a - F^\mu_{\beta a}\o^\beta_\alpha
&= F^\mu_{\alpha aj} \th^j\\ 
dF^\mu_{\alpha p} + F^\mu_{\alpha p}\o^0_0 + F^\nu_{\alpha
p}\o^\mu_\nu - F^\mu_{\alpha q}\o^q_p - F^\mu_{\beta p}\o^\beta_\alpha
&= F^\mu_{\alpha pj}\th^j
\endaligned
\end{equation}
The $F$'s satisfy the six algebraic equations GD(3.36), while their
covariant derivatives satisfy the equations GD(3.37) through GD(3.41).

The covariant derivatives of the functions in~\eqref{GD(3.27)} are
defined in a way analogous to those of the $F$'s in~\eqref{GD(3.34)},
except that the coefficient of $\omega^0_0$ must be multiplied by two
in every case.  For example,
\begin{equation}\label{GD(3.43)}
dD_{ab} + 2D_{ab}\o^0_0 - D_{cb}\o^c_a - D_{ac}\o^c_b = D_{abi}\th^i
\end{equation}
defines the covariant derivatives $D_{abi}$ of $D_{ab}$.  The
covariant derivatives of the other tensors in~\eqref{GD(3.27)} are
defined similarly.  The formula
\begin{equation}\label{GD(3.44)}
dD_{a\alpha} + 2D_{a\alpha}\o^0_0 - D_{b\alpha}\o^b_a -
D_{a\beta}\o^\beta_\alpha = D_{a\alpha i}\th^i 
\end{equation}
defines the covariant derivatives $D_{a\alpha i}$ of $D_{a\alpha}$.
The covariant derivatives of the other tensors in (3.28) are defined
similarly. 
One last set of functions, $R_i$, are defined by
\begin{equation}\label{GD(3.45)}
\o^0_{n+1} = R_i\th^i
\end{equation}
All these covariant derivatives and the $R_i$ are related in the
set of equations GD(3.46) through GD(3.59) 
in~\cite[pages 25-30]{CJ2}.  They are all real analytic functions.


\section{A sufficient condition to be isoparametric}\label{Isopara}
Consider a Legendre map $\lambda: M^{n-1}
\to \Lambda^{2n-1}$ that is proper Dupin with four distinct curvature
spheres, constant Lie curvature $r$,
and 
$M$ is a connected real analytic manifold.  For the rest of the paper
we do not use the Einstein summation convention.

By the work of M\"unzner \cite{Mun1}, \cite{Mun2}, in order for
$\lambda$ to be Lie equivalent to an isoparametric hypersurface with
four principal curvatures, it is necessary that the multiplicities of
the four curvature spheres satisfy $m_1=m_2$ and $m_3=m_4$ and that
the Lie curvature $r = -1$.  In this section, we assume these
necessary conditions and then find sufficient conditions in
Theorem~\ref{Th11} for Lie equivalence to an isoparametric
hypersurface. 

\begin{lemma}\label{L1} If the Lie curvature $r =-1$, and the
multiplicities 
satisfy $m_1=m_2$ and $m_3=m_4$, and if $Y:U \to G$ is any
second order 
Lie frame field, then on $U$ the symmetric matrices
of~\eqref{GD(3.27)} satisfy
\begin{equation}\label{scalar1}
\aligned
D_1 &= d_1 I_{m_1}\\  E_2 &= e_2 I_{m_1} \\
D_3+E_3 &= (d_3+e_3)I_{m_3}\\  D_4 - E_4 &= (d_4-e_4)I_{m_3}
\endaligned
\end{equation}
where $d_1,\dots, e_4:U \to \R$ are the real analytic functions
\begin{equation}\label{trace}
d_1 = \frac1{m_1} \trace D_1, \dots, e_4 = \frac1{m_3} \trace E_4
\end{equation}
\end{lemma}

\begin{proof} This follows from GD(3.42).
\end{proof}

For a second order Lie frame field
$Y:U\to G$, let
\begin{equation}\label{1}
\aligned
v_{a\alpha} &= (F^\alpha_{p a}, F^\mu_{\alpha a}), 
\quad |v_{a\alpha}|^2 = 2\sum_p(F^\alpha_{pa})^2 + 4
\sum_\mu(F^\mu_{\alpha a})^2 \\
v_{p\alpha} &= (F^\alpha_{p a}, F^\mu_{\alpha
p}), \quad |v_{p\alpha}|^2 = 2\sum_a(F^\alpha_{pa})^2 + 4 \sum_\mu
(F^\mu_{\alpha p})^2 \\
v_{a\mu} &= (F^\mu_{p a}, F^\mu_{\alpha a}), \quad
|v_{a\mu}|^2 = 2 \sum_p (F^\mu_{pa})^2 + 4\sum_\alpha (F^\mu_{\alpha
a})^2 \\
v_{p\mu} &= (F^\mu_{p a}, F^\mu_{\alpha p}), \quad
|v_{p\mu}|^2 
= 2\sum_a (F^\mu_{pa})^2 + 4 \sum_\alpha
(F^\mu_{\alpha p})^2
\endaligned
\end{equation}

If the Lie curvature $r =-1$, and the
multiplicities 
satisfy $m_1=m_2$ and $m_3=m_4$, then the middle four equations in
GD(3.36) 
become, when the nonsummed indices of each range are set equal,
\begin{equation}\label{3.36A1}
\aligned
|v_{a\alpha}|^2 &= d_1 - E_{aa} + E_{\alpha \alpha} \\
|v_{p\alpha}|^2 &= e_2 - D_{pp} - D_{\alpha \alpha} \\
|v_{a\mu}|^2 &= -d_1 - E_{aa} - E_{\mu\mu} \\
|v_{p\mu}|^2 &= -e_2 -D_{pp} - D_{\mu\mu}
\endaligned
\end{equation}
In addition, if all eight of the matrices $D_1,\dots, E_4$ are scalar
at every point of $U$, then equations~\eqref{3.36A1} become
\begin{equation}\label{3.36scalar}
\aligned
|v_{a\alpha}|^2 &= d_1-e_1+e_3 \\
|v_{p\alpha}|^2 &= e_2-d_2-d_3 \\
|v_{a\mu}|^2 &= -e_4-d_1-e_1 \\
|v_{p\mu}|^2 &= -e_2-d_2-d_4
\endaligned
\end{equation}
which shows that the functions on the left hand side do not depend on 
$a$, $p$, $\alpha$, or $\mu$ in this case.

\begin{remark}\label{Rem1}  If
$D_1$ is 
scalar at every point of $U$,  then a frame change~\eqref{GD(3.14)} of
the form  
\begin{equation}\label{FC}
\tilde Y = Y a(I_2, I,0,sL)
\end{equation}
can be made so that $d_1=0$ at every point of $U$.
This follows from GD(3.32), which shows that $\tilde d_1 = d_1 -s$.
\end{remark}

\begin{remark}\label{Rem2} If $D_1$ is scalar on $U$, then $D_{ab} =
d_1 
\delta_{ab}$, for all $a$ and $b$.  If we define the covariant
derivative of $d_1$ to be 
\begin{equation}\label{covder}
dd_1 + 2d_1 \omega^0_0 = \sum_i {d_1}_i \theta^i
\end{equation}
then by~GD(3.43)
\begin{equation}
\aligned
\sum_j D_{abj} \theta^j &= dD_{ab} + 2D_{ab} \omega^0_0 - \sum_c
D_{cb} \omega^c_a - \sum_c D_{ac} \omega^c_b \\
&= \delta_{ab}(dd_1 + 2d_1 \omega^0_0) - d_1(\omega^b_a + \omega^a_b)
\\
&= \delta_{ab} \sum_i d_{1i} \theta^i
\endaligned
\end{equation}
In particular, $d_{1i} = D_{aai}$, for all $a$ and $i$.  This same
principle applies to all eight of the functions $d_1,\dots,e_4$, when
all eight of these matrices are scalar.
\end{remark}

\begin{lemma}\label{LA1A2} Suppose the Lie curvature $r =-1$, and the
multiplicities 
satisfy $m_1=m_2$ and $m_3=m_4$.  If $Y:U \to G$ is
a second 
order Lie frame field for which
\begin{equation}\label{A2i}
|v_{a\alpha}| = |v_{p\mu}|, \quad \text{ and } \quad
|v_{a\mu}|=|v_{p\alpha}|
\end{equation}
for all $a$, $p$, $\alpha$,
and $\mu$,
then the eight matrices $D_1,\dots,E_4$ are scalar on
$U$ and $Y$ can be adjusted by a change~\eqref{FC} on $U$ so that
\begin{equation}\label{A2a}
d_1=0
\end{equation}
on $U$, and then
\begin{equation}\label{dande}
d_2=e_1, \; d_4 = d_3, \; e_2 = 0, \; e_3 = -d_3, \;
e_4 = d_3
\end{equation}
\begin{equation}\label{A3}
|v_{a\alpha}| = |v_{p\mu}| = |v_{p\alpha}| = |v_{a\mu}|
\end{equation}
on $U$, for all $a$, $p$, $\alpha$, and $\mu$, and
\begin{equation}\label{13} 
\omega^0_{n+1}=0
\end{equation}
on $U$.  That is, by ~\eqref{GD(3.45)}, $R_i=0$ on $U$, for all $i$.
\end{lemma}

\begin{proof} As described in Remark~\ref{Rem1}, a frame
change~\eqref{FC} will give~\eqref{A2a}, and then~\eqref{dande}
follows from GD(3.42) by linear algebra.
Putting~\eqref{dande}
into~\eqref{3.36scalar}, we obtain~\eqref{A3}.  Finally, to
prove~\eqref{13}, use GD(3.46i), for any $c$ and any $a=b$, to get
\begin{equation}\label{P1}
\aligned
d_{1c} &= D_{aac} = -R_c - \frac 2{m_1} \sum_{p,\alpha}
D_{p\alpha}(F^\alpha_{pc} + 2\delta_{ac} F^\alpha_{pa}) \\  
&-\frac2{m_1} \sum_{p,\mu} E_{p\mu}(F^\mu_{pc} +
2\delta_{ac}F^\mu_{pa}) \\ 
&+\frac8{m_1} \sum_{\alpha,p,\mu}
(F^\mu_{pa}F^\alpha_{pa}F^\mu_{\alpha c} 
+ F^\mu_{pa} F^\alpha_{pc} F^\mu_{\alpha a} + F^\mu_{pc} F^\alpha_{pa}
F^\mu_{\alpha a})
\endaligned
\end{equation}
and GD(3.62) with $a=b$ and any $c$, to get
\begin{equation}\label{P2}
\aligned
\sum_{p,\alpha}
&D_{p\alpha}(F^\alpha_{pc} + 2\delta_{ac} F^\alpha_{pa}) +
\sum_{p,\mu} E_{p\mu}(F^\mu_{pc} + 2\delta_{ac}F^\mu_{pa}) \\
&= 4 \sum_{\alpha,p,\mu} (F^\mu_{pa}F^\alpha_{pa}F^\mu_{\alpha c}
+ F^\mu_{pa} F^\alpha_{pc} F^\mu_{\alpha a} + F^\mu_{pc} F^\alpha_{pa}
F^\mu_{\alpha a})
\endaligned
\end{equation}
Substitute~\eqref{P2} into~\eqref{P1} to get
\begin{equation}\label{P3}
d_{1c} = -R_c
\end{equation}
Since $d_1=0$ on $U$, we have $d_{1j}=0$ on $U$, for every $j$,
by~\eqref{covder}.  Therefore,
\begin{equation}\label{P4}
R_c=0
\end{equation}
on $U$, for all $c$.

For any $p=q$ and for all $c$ in~GD(3.51i)
\begin{equation}\label{P5}
e_{2c} = E_{ppc} = R_c + 2\sum_\mu E_{p\mu} F^\mu_{pc} +2 \sum_\alpha
D_{p\alpha} F^\alpha_{pc}
\end{equation}
By~GD(3.51ii), we have for any $p=q$ and any $s$
\begin{equation}\label{P7}
\aligned
e_{2s} &= E_{pps} = R_s + \frac 2{m_1} \sum_{a, \alpha}
D_{a\alpha}(F^\alpha_{sa} 
+2\delta_{ps} F^\alpha_{pa}) \\ &+ \frac2{m_1} \sum_{a,\mu} D_{\mu
a}(F^\mu_{sa} + 2\delta_{ps} F^\mu_{pa}) \\
 &- \frac8{m_1} \sum_{a,\alpha,\mu} (F^\mu_{\alpha s}F^\mu_{pa}
F^\alpha_{pa} + F^\mu_{\alpha p} F^\mu_{sa} F^\alpha_{pa} +
F^\mu_{\alpha p}F^\mu_{pa} F^\alpha_{sa})
\endaligned
\end{equation}
By GD(3.63), for any $p=q$ and any $s$,
\begin{equation}\label{P8}
\aligned
&\sum_{a,\alpha} D_{a\alpha}(F^\alpha_{sa} +2\delta_{ps}F^\alpha_{pa})
+ \sum_{a,\mu} D_{\mu a}(F^\mu_{sa} + 2\delta_{ps} F^\mu_{pa}) \\
&= 4\sum_{a,\alpha,\mu} (F^\mu_{\alpha p} F^\alpha_{pa} F^\mu_{sa} +
F^\mu_{\alpha p} F^\alpha_{sa} F^\mu_{pa}+ F^\mu_{\alpha s}
F^\alpha_{pa} F^\mu_{pa})
\endaligned
\end{equation}
Substitute~\eqref{P8} into~\eqref{P7} to get
\begin{equation}\label{P9}
e_{2s} = R_s
\end{equation}
for all $s$.  Since $e_2=0$ on $U$, we have $e_{2j}=0$ on $U$, for all
$j$, and therefore
\begin{equation}\label{P10}
R_s=0
\end{equation}
for all $s$.  
By~GD(3.48iii), for any $\alpha = \beta$ and for any $\gamma$
\begin{equation}\label{P13}
\aligned
d_{3\gamma} &= D_{\alpha\alpha\gamma} = R_\gamma + \frac2{m_1}
\sum_{a,p} 
D_{pa}(F^\gamma_{pa} 
+2\delta_{\alpha\gamma} F^\alpha_{pa}) \\ &+ \frac4{m_1} \sum_{p,\mu}
E_{p\mu}(F^\mu_{\gamma p} +2 \delta_{\alpha\gamma} F^\mu_{\alpha p})
\\ &+ \frac8{m_1} \sum_{a,p,\mu} (F^\mu_{\gamma a} F^\alpha_{pa}
F^\mu_{\alpha p} + F^\mu_{\gamma p} F^\alpha_{pa} F^\mu_{\alpha a} +
F^\gamma_{pa} F^\mu_{\alpha a} F^\mu_{\alpha p})
\endaligned
\end{equation}
By GD(3.64) for all $\alpha = \beta$ and for all $\gamma$
\begin{equation}\label{P14}
\aligned
&-\sum_{a,\mu} D_{\mu a}(F^\mu_{\gamma a} + 2 \delta_{\alpha \gamma}
F^\mu_{\alpha a}) - \sum_{p,\mu}
E_{p\mu}(F^\mu_{\gamma p} +2 \delta_{\alpha\gamma} F^\mu_{\alpha p})
\\ &= 4\sum_{a,p,\mu} (F^\mu_{\gamma a} F^\alpha_{pa}
F^\mu_{\alpha p} + F^\mu_{\gamma p} F^\alpha_{pa} F^\mu_{\alpha a} +
F^\gamma_{pa} F^\mu_{\alpha a} F^\mu_{\alpha p})
\endaligned
\end{equation}
By GD(3.52iii), for any $\alpha = \beta$ and for all $\gamma$
\begin{equation}\label{P16}
\aligned
e_{3\gamma} &= E_{\alpha\alpha\gamma} = R_\gamma -\frac 2{m_1}
\sum_{a,p} 
D_{pa}(F^\gamma_{pa} 
+2\delta_{\alpha\gamma} F^\alpha_{pa}) \\&+ \frac4{m_1} \sum_{a,\mu}
D_{\mu 
a}(F^\mu_{\gamma a} + 2 \delta_{\alpha \gamma} F^\mu_{\alpha a}) \\
&+ \frac8{m_1} \sum_{a,p,\mu} (F^\mu_{\gamma a} F^\alpha_{pa}
F^\mu_{\alpha p} + F^\mu_{\gamma p} F^\alpha_{pa} F^\mu_{\alpha a} +
F^\gamma_{pa} F^\mu_{\alpha a} F^\mu_{\alpha p})
\endaligned
\end{equation}
Now $e_3+d_3=0$ on $U$, so $e_{3\gamma} + d_{3\gamma}=0$ on $U$, for
all $\gamma$, so~\eqref{P13} and~\eqref{P16} added together give
on $U$, for all $\gamma$,
\begin{equation}\label{P17a}
\aligned
0 &= 2R_\gamma \\ &+ \frac4{m_1} \left( \sum_{p,\mu}
E_{p\mu}(F^\mu_{\gamma p} +2 \delta_{\alpha\gamma} F^\mu_{\alpha p})
+ \sum_{a,\mu} D_{\mu a}
(F^\mu_{\gamma a} + 2 \delta_{\alpha \gamma} F^\mu_{\alpha a})\right)
\\ 
&+ \frac{16}{m_1} \sum_{a,p,\mu} (F^\mu_{\gamma a} F^\alpha_{pa}
F^\mu_{\alpha p} + F^\mu_{\gamma p} F^\alpha_{pa} F^\mu_{\alpha a} +
F^\gamma_{pa} F^\mu_{\alpha a} F^\mu_{\alpha p})
\endaligned
\end{equation}
and this with~\eqref{P14} implies
\begin{equation}\label{P18}
R_\gamma = 0
\end{equation}
on $U$ for all $\gamma$.
In the same way, by~GD(3.48iv)
\begin{equation}\label{P18a}
d_{3\mu} = D_{\alpha\alpha\mu} = R_\mu - D_{\mu\alpha\alpha}+ 6\sum_a
D_{a\alpha}F^\mu_{\alpha a} + 2\sum_p D_{p\alpha}F^\mu_{\alpha p}
\end{equation}
By GD(3.52iv),
\begin{equation}\label{P18b}
e_{3\mu} = E_{\alpha\alpha\mu} = R_\mu +D_{\mu\alpha\alpha} + 2\sum_a
D_{a\alpha}F^\mu_{\alpha a} + 6\sum_p D_{p\alpha} F^\mu_{\alpha p}
\end{equation}
Adding these equations together and using~\eqref{P18} and the fact
that $d_{3\gamma}+e_{3\gamma}=0$ on $U$,
we get on $U$, for every $\alpha$ and $\mu$,
\begin{equation}\label{P18c}
\sum_a D_{a\alpha} F^\mu_{\alpha a} + \sum_p D_{p\alpha} F^\mu_{\alpha
p} = 0
\end{equation}
Finally, by~GD(3.49iv), for any
$\mu=\nu$ and for any $\sigma$,
\begin{equation}\label{P19}
\aligned
d_{4\sigma} &= D_{\mu\mu\sigma} = -R_\sigma + \frac2{m_1} \sum_{a,p}
D_{pa}(F^\sigma_{pa} 
+2\delta_{\mu\sigma} F^\mu_{pa}) \\ &+ \frac4{m_1} \sum_{p,\alpha}
D_{p\alpha} 
(F^\sigma_{\alpha p} + 2 \delta_{\mu\sigma}F^\mu_{\alpha p})\\
&+\frac8{m_1} \sum_{a,p,\alpha} (F^\sigma_{\alpha a}
F^\mu_{pa}F^\mu_{\alpha p} + F^\sigma_{\alpha p} F^\mu_{pa}
F^\mu_{\alpha a} + F^\sigma_{pa}F^\mu_{\alpha a} F^\mu_{\alpha p})
\endaligned
\end{equation}
By GD(3.65) with $\mu=\nu$ and for any $\sigma$,
\begin{equation}\label{P20}
\aligned
&\sum_{a,\alpha} D_{a\alpha}(F^\sigma_{\alpha a} +2\delta_{\mu\sigma}
F^\mu_{\alpha a}) + \sum_{p,\alpha} D_{p\alpha}
(F^\sigma_{\alpha p} + 2 \delta_{\mu\sigma}F^\mu_{\alpha p}) \\
&= -4 \sum_{a,p,\alpha} (F^\sigma_{\alpha a}
F^\mu_{pa}F^\mu_{\alpha p} + F^\sigma_{\alpha p} F^\mu_{pa}
F^\mu_{\alpha a} + F^\sigma_{pa}F^\mu_{\alpha a} F^\mu_{\alpha p})
\endaligned
\end{equation}
By GD(3.53iv), for all $\mu=\nu$ and for any $\sigma$,
\begin{equation}\label{P22}
\aligned
e_{4\sigma} &= E_{\mu\mu\sigma}= R_\sigma + \frac2{m_1} \sum_{a,p}
D_{pa}(F^\sigma_{pa} 
+2\delta_{\mu\sigma} F^\mu_{pa}) \\ &- \frac4{m_1} \sum_{a,\alpha}
D_{a\alpha}(F^\sigma_{\alpha a} +2\delta_{\mu\sigma} 
F^\mu_{\alpha a}) \\
&-\frac 8{m_1} \sum_{a,p,\alpha} (F^\sigma_{\alpha a}
F^\mu_{pa}F^\mu_{\alpha p} + F^\sigma_{\alpha p} F^\mu_{pa}
F^\mu_{\alpha a} + F^\sigma_{pa}F^\mu_{\alpha a} F^\mu_{\alpha p})
\endaligned
\end{equation}
Now $e_4 = d_4$ on $U$ implies
that $e_{4\sigma}-d_{4\sigma}=0$ on $U$, so by~\eqref{P19}
and~\eqref{P22} we get
\begin{equation}\label{P23a}
\aligned
0 &= e_{4\sigma} - d_{4\sigma} 
= 2R_\sigma \\ &- \frac 4{m_1}
\left(\sum_{a,\alpha} D_{a\alpha}(F^\sigma_{\alpha a} +
2\delta_{\mu\sigma} F^\mu_{\alpha a}) + \sum_{p,\alpha}
D_{p\alpha} (F^\sigma_{\alpha p} + \delta_{\mu\sigma} F^\mu_{\alpha
p}\right) \\ &- \frac{16}{m_1} \sum_{a,p,\alpha} (F^\sigma_{\alpha a}
F^\mu_{pa}F^\mu_{\alpha p} + F^\sigma_{\alpha p} F^\mu_{pa}
F^\mu_{\alpha a} + F^\sigma_{pa}F^\mu_{\alpha a} F^\mu_{\alpha
p}) \\ 
&= 2R_\sigma
\endaligned
\end{equation}
by~\eqref{P20}.  Therefore, on $U$, 
\begin{equation}\label{P24}
R_\sigma =0
\end{equation}
for every $\sigma$.  Therefore,~\eqref{13} holds by~\eqref{P4},
\eqref{P10}, \eqref{P18}, and~\eqref{P24}.
\end{proof}

\begin{theorem}\label{Th11} Suppose the
multiplicities satisfy
$m_1=m_2$, $m_3=m_4$, and the Lie curvature is $r=-1$.  
Suppose that for any
point in $M$ there exists a second 
order frame 
field $Y:U\to G$ along $\lambda$ on an open set $U \subset M$ about
the point,  such that equations~\eqref{A2i} hold on $U$, for
all $a,p,\alpha,\mu$.  If for some $a, \a$
\begin{equation}\label{2}
|v_{a\alpha}|  > 0
\end{equation}
on an open dense subset of $U$; and if
\begin{equation}\label{2a}
d\omega^0_0 = 0
\end{equation}
on $U$, then
$\lambda:M\to \Lambda$ is Lie equivalent to an isoparametric
hypersurface.
\end{theorem}

\begin{proof}  Given any point of $M$,
let $Y:U \to G$ be a second order frame field about the point
satisfying the hypotheses.  
By Lemma~\ref{LA1A2}, we may assume $Y$ satisfies~\eqref{A2a},
\eqref{dande}, and~\eqref{A3} on $U$.  
Thus,~\eqref{2} implies that all the
functions in~\eqref{A3} are positive on $U$.
These properties are preserved
by any frame change of the form
\begin{equation}\label{X1}
\tilde Y = Y a(tI_2, I,0,0)
\end{equation}
where $t$ is any nowhere zero real analytic function on $U$, in which
case 
\begin{equation}\label{X2}
\tilde \omega^0_0 = \omega^0_0 + d\log|t|
\end{equation}
on $U$.  We may assume that $U$ is contractible.
Then~\eqref{2a} implies that
\begin{equation}\label{X3}
\omega^0_0 = df
\end{equation}
for some real analytic function $f$ on $U$.  Making the frame
change~\eqref{X1} with $t = e^{-f}$, we have
\begin{equation}\label{X4}
\tilde \omega^0_0 = 0
\end{equation}
on $U$.  We now continue with this frame and drop the tildes.
By~\eqref{GD(3.33)}, our hypothesis $d\omega^0_0=0$ on $U$ implies that
$D_{a\alpha}$ and its covariant derivatives $D_{a\alpha j}$ are
identically zero on $U$.  Then~GD(3.54) with~\eqref{A2a}
and~\eqref{dande} implies that 
\begin{equation}\label{Co1}
\aligned
0&= (d_3-d_2)F^\alpha_{pa} \\
0 &= (d_4+e_4-2e_1-e_3-d_3)F^\mu_{\alpha a} = 2(d_3-d_2)F^\mu_{\alpha
a}
\endaligned
\end{equation}
on $U$, for any $a$, $p$, $\alpha$, and $\mu$.  Thus,
\begin{equation}\label{T2aa}
0= (d_3 - d_2)^2|v_{a\alpha}|^2
\end{equation}
at every point of $U$, for all $a, \alpha$, and so~\eqref{2} implies
\begin{equation}\label{T2a}
d_3-d_2 =0
\end{equation}
on an open dense subset of $U$, hence on all of $U$, by continuity. 
So,~\eqref{dande} becomes
\begin{equation}\label{TTT2b}
d_1=e_2=0, \quad e_1=d_2=d_3=d_4=e_4=-e_3
\end{equation}
on $U$.  Since $d\omega^0_0 = 0$ and $\omega^0_{n+1}=0$ on $U$, we get
from GD(3.47) that 
\begin{equation}\label{X5}
d_{2a} = d_{2\alpha} = d_{2\mu}=0
\end{equation}
on $U$, for all $a,\alpha,\mu$, and from GD(3.48) that
\begin{equation}\label{X6}
d_{3p}=0
\end{equation}
on $U$, for all $p$.  Since $d_2 = d_3$ on $U$, it follows that $d_2$
is covariant constant on $U$, and therefore $d_2$ must be constant on
$U$, since $\omega^0_0 = 0$ in~\eqref{covder} and $U$ is connected.
Putting~\eqref{TTT2b} into~\eqref{3.36scalar} and using~\eqref{2}, we
have 
\begin{equation}\label{TT2b}
-2d_2 = |v_{a\alpha}|^2 >0
\end{equation}
on an open dense subset of $U$.
Therefore, $d_2$ is a negative constant.
Making a frame change~\eqref{X1} with the constant
\begin{equation}\label{TT2c}
t = \sqrt{-d_2}
\end{equation}
we have, by~GD(3.32), that
\begin{equation}\label{TT2d}
\tilde d_2 = \frac1{t^2} d_2 = -1
\end{equation}
at every point of $U$.  Hence, we may assume that
\begin{equation}\label{T2}
d_2=-1
\end{equation}
on $U$.  We have thus proved that for any point of $M$, there exists a
second order frame field $Y:U\to G$ on a neighborhood of that point
for which $\omega^0_0 = 0$, $\omega^0_{n+1}=0$, and
\begin{equation}\label{T1}
d_1=e_2=0, \quad e_1=d_2=d_3 = d_4 = e_4 = -e_3 = -1
\end{equation}
on $U$.
The following equations then follow from the structure
equations 
\[
dY_A = \sum_{B=0}^{n+2} \omega^B_A Y_B
\]
and the properties of our frame.
\begin{equation}\label{T3}
\aligned
dY_0 &= \sum_p \theta^p Y_p + \sum_\alpha \theta^\alpha Y_\alpha +
\sum_\mu \theta^\mu Y_\mu \\
dY_1 &= \sum_a \theta^a Y_a - \sum_\alpha \theta^\alpha Y_\alpha +
\sum_\mu \theta^\mu Y_\mu \\
dY_{n+1} &= -\sum_a \theta^a Y_a + \sum_\alpha \theta^\alpha Y_\alpha
- \sum_\mu \theta^\mu Y_\mu \\
dY_{n+2} &= -\sum_p \theta^p Y_p - \sum_\alpha \theta^\alpha Y_\alpha
- \sum_\mu \theta^\mu Y_\mu
\endaligned
\end{equation}
If we let
\begin{equation}\label{T4}
W_1 = Y_0 + Y_{n+2}, \quad W_2 = Y_1 + Y_{n+1}
\end{equation}
then equations~\eqref{T3} show that
\begin{equation}\label{T5}
dW_1=0, \quad dW_2 = 0
\end{equation}
on $U$, so $W_1$ and $W_2$ are constant vectors (assuming $U$
connected).  In addition, they span a time-like line in $\R^{n+3}_2$,
since
\begin{equation}\label{T6}
\aligned
\langle W_1, W_1 \rangle &= -2 \\
\langle W_2, W_2 \rangle &= -2 \\
\langle W_1, W_2 \rangle &= 0
\endaligned
\end{equation}
Then $W_1$, $W_2$, $W_1-W_2$, and $W_1+W_2$ are four points on this
time-like line such that
\begin{equation}\label{T7}
\aligned
\langle Y_0, W_1 \rangle &= 0 \\
\langle Y_1, W_2 \rangle &=0 \\
\langle Y_0 + Y_1, W_1 - W_2 \rangle &= 0 \\
\langle -Y_0 + Y_1, W_1+W_2 \rangle &=0
\endaligned
\end{equation}
on $U$.  
If $\tilde Y:\tilde U \to G$ is another Lie frame field
satisfying~\eqref{T1}, then on the intersection $\tilde
U \cap U$ (supposed nonempty) they must be related by
\[
\tilde Y = Y a(I_2, B,0,0)
\]
where $B:U\cap \tilde U \to O(n-1)$ is an analytic map.  In particular,
\[
\tilde Y_0 = Y_0,\; \tilde Y_1 = Y_1, \; \tilde Y_{n+1} = Y_{n+1}, \;
\tilde Y_{n+2} = Y_{n+2}
\]
and therefore~\eqref{T7} holds for $\tilde Y$, for the same vectors
$W_i$, for $i = 1,2,3,4$, and thus~\eqref{T7} holds on all of $M$ for
the four curvature spheres.
By Cecil's Theorem~5.6 (\cite[pp
102-103]{Cec6}), 
$\lambda:M\to \Lambda$ is Lie equivalent to the Legendre submanifold
induced by an isoparametric hypersurface.  
\end{proof}

\section{Reducibility}\label{Necessary}
Before we return to the case of a Dupin hypersurface with four
principal curvatures, we prove some general results about reducible
Dupin hypersurfaces.

Pinkall \cite{Pinkall} introduced the basic constructions of building
tubes, cylinders, cones, and surfaces of revolution over a Dupin
hypersurface $M^{n-1}$ in $\R^n$ with $g$ principal curvatures to get
a Dupin hypersurface $W^{n-1-k}$ in $\R^{n+k}$ with $g+1$ principal
curvatures.  In general, these constructions introduce a new principal
curvature of multiplicity $k$, which is easily seen to be constant
along its curvature surfaces.  The other principal curvatures are
determined by the principal curvatures of $M^{n-1}$, and the Dupin
property is preserved for these principal curvatures.  A Dupin
hypersurface that is locally Lie equivalent to a hypersurface
obtained by one of these constructions is said to be reducible.  In
Theorem~4 of his paper, Pinkall gave a formulation of reducibility in
terms of Lie sphere geometry.  As in the paper \cite{CJ}, we use this
formulation as our definition of reducibility on an open subset of a
Dupin submanifold as follows.

\begin{definition}\label{reducible}
We define the
Dupin submanifold $\lambda:M \to \Lambda$ to be \textit{reducible on an
open subset} $O\subset M$ if
some curvature
sphere maps $O$ into some fixed linear subspace
of $\mathbf{RP}^{n+2}$ of codimension at least two.  We say that
$\lambda$ is \textit{reducible} if it is reducible on $M$.
Define $\lambda$
to be \textit{locally irreducible} if it is not reducible on any open
subset of $M$.  Define $\lambda$ to be \textit{irreducible} if it is
not reducible.
\end{definition}

\begin{proposition}\label{5P2}  If a connected,
proper Dupin submanifold is reducible on an open subset $U$ of $M$,
then it is reducible.  Thus, a connected, proper Dupin submanifold is
locally irreducible if and only if it is 
irreducible. 
\end{proposition}

\begin{proof} Let $\lambda: M^{n-1} \to \Lambda^{2n-1}$ be a
connected, proper Dupin submanifold.  Suppose that there exists an
open subset $U \subset M$ such that the restriction of $\lambda$ to
$U$ is reducible.  By the definition of reducibility, there exists a
curvature sphere $[K]$ of $\lambda$ and two linearly independent
vectors $v_i \in \R^{n+3}_2$, $i=1,2$, such that
\begin{equation*}
\langle K,v_i \rangle =0
\end{equation*}
on the set $U$.  Since $\lambda$ is analytic, the curvature sphere map
$[K]$ is analytic on $M$, and so the functions $\langle K, v_i
\rangle$ are analytic on $M$.  Since these functions equal zero on the
open set $U$, they are equal to zero on all of the connected manifold
$M$, and thus $\lambda: M \to \Lambda$ is reducible.
\end{proof}

This result has ramifications for the case of proper Dupin
hypersurfaces with $g=3$ principal curvatures.

\begin{corollary}\label{5C1} Let $\lambda: M^{n-1} \to \Lambda^{2n-1}$
be an irreducible proper Dupin submanifold with $g=3$ principal
curvatures.  Then $\lambda$ is Lie equivalent to an isoparametric
hypersurface in $S^n$.
\end{corollary}

\begin{proof} In the case where all the principal curvatures have
multiplicity one, this was proven by Pinkall \cite{Pink}.  For the
case of higher multiplicities, this was proven in \cite[p.\ 175]{CJ}
under the assumption that $\lambda$ is locally irreducible.  By
Proposition~\ref{5P2}, we see that the hypothesis of irreducibility is
sufficient.
\end{proof}

We now prove a characterization of reducibility for proper
Dupin submanifolds.  Note that we need only use three of Pinkall's
four constructions, since as Pinkall showed, the cone construction is
locally Lie equivalent to the tube construction.

\begin{proposition}\label{5P1} Let $\nu:W^{d-1} \to \Lambda^{2d-1}$ be
a connected,  reducible proper Dupin submanfold.  Then $\nu$
is Lie equivalent to a proper Dupin submanifold $\mu$ which is
obtained from a lower dimensional proper Dupin submanifold $\lambda$
by one of Pinkall's three constructions (tube, cylinder, surface of
revolution).
\end{proposition}

\begin{proof} It is possible that the curvature sphere $[K]$ of $\nu$
locally 
lies in a linear subspace of codimension even higher than two.  For
each $x \in W$, let $m_x$ be the largest positive integer such that
for some neighborhood $U_x$, the curvature sphere map $[K]$ restricted
to $U_x$ is contained in a linear space of codimension $m_x+1$ in
$\mathbf{RP}^{d+2}$. By hypothesis, we know that $m_x \geq 1$ for all
$x \in W$.  Choose $x_0$ to be a point where $m_x$ attains its maximum
value $m$.  Then there exist linearly independent vectors
$v_1,\dots,v_{m+1}$ in $\R^{d+3}_2$ such that on an open set $U_{x_0}$
about $x_0$,
\begin{equation}\label{5e1}
\langle K,v_i \rangle = 0
\end{equation}
for $1 \leq i \leq m+1$.  Since $\nu$ is analytic, the curvature
sphere map $[K]:W \to Q$ is analytic, and since the analytic functions
$\langle K, v_i \rangle$ equal zero on the open set $U_{x_0}$, they
must equal zero on the whole connected manifold $W$.
Thus,~\eqref{5e1} holds on all of $W$, and the function $m_x = m$ for
all $x \in W$.

The rest of the proof is essentially the same as the proof of
Theorem~2.8 in \cite[pp.\ 145-147]{Cec6}.  Specifically, let $E$ be the
linear subspace in $\R^{d+3}_2$ of codimension $m+1$ whose orthogonal
complement $E^\perp$ is spanned by the vectors $v_1,\dots,v_{m+1}$.
The signature of $\langle\, ,\, \rangle$ on $E^\perp$ must be
$(m+1,0)$, $(m,1)$, or $(m,0)$.  Then, as in the proof of Theorem~2.8,
one can show that there is a Lie transformation $A$ such that $\mu =
A\nu$ is obtained from a proper Dupin submanifold $\lambda:M^{n-1} \to
\Lambda^{2n-1}$, where $n = d-m$, by the surface of revolution, tube,
or cylinder construction, depending on whether the signature of the
inner product on $E^\perp$ is $(m+1,0)$, $(m,1)$, or $(m,0)$,
respectively.  The proof of Theorem~2.8 deals specifically with the
case where $[K]$ has multiplicity $m$, and so $\mu$ has one more
curvature sphere than $\lambda$.  In the case where the multiplicity
of $[K]$ is greater than $m$, one must make some slight adjustments in
the exposition of the proof.  In that case, the curvature sphere
$A[K]$ is equal to one of the curvature spheres of $\mu$ induced from
a curvature sphere $[k]$ of $\lambda$, and the multiplicity of $[K]$
is $m+q$, where $q$ is the multiplicity of $[k]$ as a curvature sphere
of $\lambda$.  In that case, $\mu$ and $\lambda$ have the same number
of distinct curvature spheres. 
\end{proof}

\begin{remark}\label{Rem22} Pinkall \cite[p.\ 438]{Pinkall} proved
that $\nu$ as in 
Proposition~\ref{5P1} is locally Lie equivalent to a proper Dupin
submanifold $\mu$ that is obtained by one of the three constructions.
By using analyticity, we are able to eliminate the word
\textit{locally} from the statement of the result.
\end{remark}


We return to the case where $\lambda: M^{n-1} \to \Lambda^{2n-1}$
is a proper Dupin submanifold with four curvature spheres of
multiplicities $m_1=m_2$, $m_3 = m_4$, and with constant Lie curvature
$r=-1$.

\begin{definition}\label{DS1} For a second order frame field $Y:U \to
G$, let $\1'$ be the set of all $a \in \1$ such that
\begin{equation}\label{S1}
F^\alpha_{pa}=0, \; F^\mu_{pa} = 0,\; F^\mu_{\alpha a} = 0
\end{equation}
on $U$, for all $p,\alpha,\mu$.
\end{definition}

\begin{proposition}\label{PS1}  If $m_1 = m_2$, $m_3=m_4$, and $r=-1$,
and if $\1'$ is nonempty for the second order frame field $Y:U
\to G$, then the symmetric matrices $D_1$, $E_2$, $D_3$, $D_4$, $E_3$,
and 
$E_4$ are scalar matrices,
\begin{equation}\label{S1a}
E_1 = \begin{pmatrix} (d_1+e_3)I_m & 0 \\ 0 & * \end{pmatrix}
\end{equation}
where $m$ is the cardinality of $\1'$, and
\begin{equation}\label{S1b}
d_1+e_2 = 0
\end{equation}
at every point of $U$.  Here $d_1,\dots,e_4$ are defined
in~\eqref{trace}. 
\end{proposition}

\begin{proof} By Lemma~\ref{L1}, just the assumptions on the
multiplicities and on $r$ imply that $D_1$, $E_2$, $D_3+E_3$, and
$D_4-E_4$ are scalar matrices on $U$.  They are given
by~\eqref{scalar1}.  For each $a \in \1'$ and $e \in \1$, the left side
of~GD(3.36ii) is zero, so that
\begin{equation}\label{S2}
0 = 2\sum_p F^\alpha_{pa} F^\beta_{pa} + 4\sum_\mu F^\mu_{\alpha a}
F^\mu_{\beta a} = (d_1 - E_{aa})\delta_{\alpha \beta} + E_{\alpha \beta}
\end{equation}
on $U$ for all $\alpha, \beta$.  Therefore, $E_3= e_3 I_{m_3}$ is a
scalar 
matrix, where
\begin{equation}\label{S3}
e_3  = E_{aa} - d_1
\end{equation}
for every $a \in \1'$.  Since $D_3+E_3$ is scalar, it follows that
\begin{equation}\label{S4}
D_3 = d_3 I_{m_3}
\end{equation}
is scalar also.  Then for any $a\in \1'$, $e \in \1$, and $\alpha =
\beta$ in~GD(3.36ii), we have
\begin{equation}\label{S4a}
0 = 2\sum_p F^\alpha_{pa}F^\alpha_{pe} + 4\sum_\mu F^\mu_{\alpha a}
F^\mu_{\alpha e} = (d_1 - E_{ae}) + e_3 \delta_{ae}
\end{equation}
on $U$.  Therefore,
\begin{equation}\label{S4b}
E_{ae} = d_1 + e_3 \delta_{ae}
\end{equation}
from which~\eqref{S1a} follows.
In the same way, for all $a=b \in \1'$, the left side
of~GD(3.36iv) is zero, and so we have
\begin{equation}\label{S5}
0 = -2\sum_p F^\mu_{pa} F^\nu_{pa} -4 \sum_\alpha F^\nu_{\alpha a}
F^\nu_{\alpha a} = E_{\mu\nu} + (2d_1+e_3)\delta_{\mu\nu}
\end{equation}
on $U$, for all $\mu, \nu$.  Therefore, $E_4 = e_4 I_{m_3}$ is a
scalar matrix, where
\begin{equation}\label{S6}
e_4 = -(2d_1 + e_3)
\end{equation}
Since $D_4 - E_4$ is a scalar matrix, it follows that
\begin{equation}\label{S7}
D_4 = d_4 I_{m_3}
\end{equation}
is scalar also.  Finally, if $a \in \1'$, then the left side of
GD(3.36i) is zero, so~\eqref{S1b} holds.
\end{proof}

\begin{corollary}\label{CS1}
Under the hypotheses of Proposition~\ref{PS1}, the second order frame
$Y:U 
\to G$ can be chosen so that
\begin{equation}\label{S8}
d_1=0,\; e_2=0,\; e_3 = -e_4,\; d_3=d_4,\; e_a = e_3,\; \forall a \in
\1' 
\end{equation}
\end{corollary}

\begin{proof}
If we make a change of frame of the form~\eqref{GD(3.14)}
\begin{equation}\label{S9}
\tilde Y = Y a(tI_2, B, 0, sL)
\end{equation}
with $t=1$, $B=I$ and $s = d_1$, then by GD(3.32) and~\eqref{S1b}
\begin{equation}\label{S10}
\tilde d_1 = 0 = \tilde e_2
\end{equation}
Dropping the tildes, we see that two of the remaining equations
in~\eqref{S8} then follow from~\eqref{S3} and~\eqref{S6}.  It remains
to prove that $d_3 = d_4$ in this frame.  For this we use GD(3.42).
In fact, using the already established equations in~\eqref{S8} and
GD(3.42ii), we have
\begin{equation}\label{S11}
(m_1+\frac{m_3}2)(d_3 + e_3) = m_1(e_1-d_2) - \frac{m_3}2 (d_4 +e_3)
\end{equation}
From GD(3.42iv), we have
\begin{equation}\label{S12}
(m_1 + \frac{m_3}2)(-e_3-d_4) = m_1(-e_1+d_2) + \frac{m_3}2(e_3+d_3)
\end{equation}
Adding \eqref{S11} and \eqref{S12}, we get
\begin{equation}\label{S13}
(m_1+\frac{m_3}2)(d_3-d_4) = \frac{m_3}2(d_3-d_4)
\end{equation}
from which we conclude that $d_3=d_4$.
\end{proof}

We shall call a second order frame $Y:U \to G$ \textit{normalized} if
it satisfies~\eqref{S8}.

For a second order frame $Y:U \to G$, let
\begin{equation}\label{3}
\aligned
f &= 2\sum_{a,p,\alpha}(F^\alpha_{pa})^2, \quad &g =
2\sum_{a,p,\mu}(F^\mu_{pa})^2\\
h &= 2\sum_{a,p,\mu}(F^\mu_{\alpha a})^2, \quad &k =
2\sum_{p,\alpha,\mu} (F^\mu_{\alpha p})^2
\endaligned
\end{equation}
By GD(3.10), a change of second order frame field~\eqref{GD(3.14)}
multiplies these
functions by a nowhere zero function.  In particular, the zero sets of
these functions are globally well defined.

\begin{theorem}\label{TS2} Suppose $\lambda:M^{n-1} \to
\Lambda^{2n-1}$ is proper Dupin with multiplicities $m_1=m_2$,
$m_3=m_4$, and Lie curvature $r = -1$.
If at least three of
$f$, $g$, $h$, and $k$ are zero 
on $M$, then $\lambda:M \to \Lambda^{2n-1}$ is reducible.
\end{theorem}

\begin{proof}
Suppose $f=g=h=0$ at every point of $M$.
Then
\begin{equation}\label{R41}
F^\alpha_{pa} = F^\mu_{pa} = F^\mu_{\alpha a} = 0
\end{equation}
for all $a,p,\alpha,\mu$ at each point of $U$
for any second order frame field $Y:U\to G$.  
Note that~\eqref{R41} is equivalent to
\begin{equation}\label{R41aa}
\omega^p_a = \omega^\alpha_a = \omega^\mu_a = 0
\end{equation}
on $U$, for all $a, p, \alpha,\mu$.
Because $\1'=\1$
in this case, we may assume that $Y$ is normalized, and
then~\eqref{S8} implies that $E_1 = e_1 I_m$ is scalar with
$e_1=e_3$. By GD(3.37), GD(3.38), and GD(3.39) together
with~\eqref{R41}, we have that
\begin{equation}\label{S14}
D_{p\alpha} = D_{a\alpha} = D_{pa} = E_{p\mu} = D_{\mu a} = D_{\mu
\alpha} =0
\end{equation}
at every point of $U$.  Therefore,
\begin{equation}\label{S15}
d \omega^0_0 = 0,\; \omega^0_a = 0 = \omega^1_p
\end{equation}
at every point of $U$.  Moreover, 
\begin{equation}\label{S16}
\omega^0_{n+1} = 0
\end{equation}
on $U$, because
GD(3.46i), GD(3.46ii),
and~\eqref{R41} imply that $R_a = R_p = 0$ on $U$, and GD(3.46iii),
GD(3.46iv), and~\eqref{R41} imply that $R_\alpha = D_{a\alpha a} = 0$
and $R_\mu = D_{\mu a a} = 0$ on $U$.  It then follows from GD(3.51)
that $e_3$ is covariant constant; that is,
\begin{equation}\label{S16a}
de_3 + 2e_3\omega^0_0 = 0
\end{equation}
on $U$.

Let $V(u)$ be the subspace of $\R^{n+3}_2$ defined by the span of the
vectors 
\begin{equation}\label{S17}
Y_0, Y_p, Y_\alpha, Y_\mu, Y_{n+2}, e_3 Y_1 - Y_{n+1}
\end{equation}
for all $p,\alpha,\mu$ at the point $u\in U$.  Then $V(u)$ does not
depend on the choice of normalized second order frame field at $u$,
since any other is given by~\eqref{GD(3.14)} with $s=0$.
Let $V$ be the span of
$V(u)$ for all $u \in M$.  
We want to prove that
$V$ is a subspace of codimension $m+1$, because the curvature sphere
$[Y_0]$ on $U$ takes all of its values in $V$ then shows that
$\lambda$ is reducible on $U$, and therefore $\lambda$ is reducible by
Proposition~\ref{reducible}.  
Since the codimension of
$V(u)$ is $m+1$ for any $u \in M$, we will obtain our result if we
prove that $V(u)$ is constant on the domain $U$ of any normalized
second order frame field $Y$.  This will be true if we show that the
derivatives of the vectors spanning $V(u)$ are zero modulo $V(u)$, for
every $u \in U$.  This follows from~\eqref{S15}, and~\eqref{S16}
and~\eqref{S16a}, and~\eqref{R41aa}.  In fact, if $\equiv$ denotes
equality modulo $V(u)$, then
\begin{equation}\label{S18}
\aligned
dY_0 &\equiv 0 \\
dY_p &\equiv \omega^1_p Y_1 + \omega^a_p Y_a + \omega^{n+1}_p Y_{n+1}
= 0 \\
dY_\alpha &\equiv \omega^1_\alpha Y_1 + \omega^a_\alpha Y_a +
\omega^{n+1}_\alpha Y_{n+1} \\
&= (e_3 Y_1 - Y_{n+1}) \theta^\alpha \equiv 0 \\
dY_\mu &\equiv \omega^1_\mu Y_1 + \omega^a_\mu Y_a + \omega^{n+1}_\mu
Y_{n+1} \\
&= (-e_3 Y_1 + Y_{n+1}) \theta^\mu \equiv 0\\
dY_{n+2} &\equiv \omega^1_{n+2} Y_1 + \omega^a_{n+2} Y_a +
\omega^{n+1}_{n+2} Y_{n+1} = 0\\
d(e_3 Y_1 - Y_{n+1}) &\equiv de_3 Y_1 + e_3(\omega^1_1 Y_1 + \theta^a
Y_a) - (\omega^a_{n+1} Y_a + \omega^{n+1}_{n+1}Y_{n+1}) \\
&= (de_3 + e_3 \omega^0_0) Y_1 + \omega^0_0 Y_{n+1} = (-e_3 Y_1 +
Y_{n+1})\omega^0_0 \equiv 0
\endaligned
\end{equation}
\end{proof}  

\begin{theorem}\label{PU1} Let $\lambda:M \to \Lambda$ be a proper
Dupin 
hypersurface for which $r=-1$, $m_1=m_2$, and $m_3=m_4$.  If $Y:U\to
G$ is a second order frame field such that
\begin{equation}\label{U4}
d\omega^0_0 = 0
\end{equation}
on $U$, and
\begin{equation}\label{U5}
\1' \neq \emptyset
\end{equation}
(see Definition~\ref{DS1}), then there exists a nonempty open subset
$W$ of $U$ such that 
the curvature sphere $[Y_0]:W \to \mathbf{RP}^{n+2}$ takes values
in a constant linear subspace of codimension at least 2.  Thus,
$\lambda$ is reducible on $W$. 
\end{theorem}

\begin{proof} By Proposition~\ref{PS1} and Corollary~\ref{CS1}, we may
assume~\eqref{S8} and~\eqref{S1a} hold.  Choosing $U$ to be
contractible, we may assume 
\begin{equation}\label{U6}
\omega^0_0 = 0
\end{equation}
on $U$.  In fact, if $U$ is contractible, then $\omega^0_0 = df$, for
some function $f$ on $U$.  A change of frame $\tilde Y = Y
a(tI_2,I,0,0)$ doesn't affect the assumptions already made, and
$\tilde \omega^0_0 = \omega^0_0 + dt/t$.  Thus, take $t = e^{-f}$.  
There exists a dense open subset of $U$ on which
$Y$ can be chosen to diagonalize $E_1$ and $D_2$.  Let $\tilde U$ be
an open connected component of this open dense subset.
It follows from GD(3.46) and~\eqref{U4}, that for all $a \in \1'$
and for all $j$,
\begin{equation}\label{U7}
0 = d_{aj} = D_{aaj}= -R_j
\end{equation}
on $\tilde U$, and hence
\begin{equation}\label{U8}
\omega^0_{n+1} = 0
\end{equation}
on $\tilde U$.  From~\eqref{U4} and~\eqref{GD(3.33)}, we know that the
six 
sets of invariants $D_{a\alpha}$, etc. and their covariant derivatives
are identically zero on $U$.  Using~\eqref{S8}, we get from GD(3.54)
\begin{equation}\label{GD354}
\aligned
0 &= D_{e\alpha p} = (-d_p + d_3) F^\alpha_{pe} \\
0 &= D_{e\alpha \mu} = -2(e_e + e_3)F^\mu_{\alpha e}
\endaligned
\end{equation}
on $\tilde U$, for all $e,p,\alpha, \mu$; from GD(3.55)
\begin{equation}\label{GD355}
\aligned
0&= D_{pe\alpha} = (e_e + e_3 +d_p - d_3) F^\alpha_{pe} \\
0 &= D_{pe \mu} = (e_e + e_3 + d_p - d_3) F^\mu_{pe}
\endaligned
\end{equation}
on $\tilde U$, for all $e,p,\alpha,\mu$; from GD(3.56)
\begin{equation}\label{GD356}
\aligned
0 &= D_{p\alpha e} = (e_e + e_3) F^\alpha_{pe} \\
0 &= D_{p\alpha \mu} = 2(d_p - d_3) F^\mu_{\alpha p}
\endaligned
\end{equation}
on $\tilde U$, for all $e,p,\alpha,\mu$; from GD(3.57)
\begin{equation}\label{GD357}
\aligned
0 &= D_{\mu e p} = (d_p - d_3) F^\mu_{pe} \\
0 &= D_{\mu e \alpha} = -2(e_e + e_3) F^\mu_{\alpha e}
\endaligned
\end{equation}
on $\tilde U$, for all $e,p,\alpha,\mu$; from GD(3.58)
\begin{equation}\label{GD358}
\aligned
0 &= D_{\mu\alpha e} = -2(e_e + e_3) F^\mu_{\alpha e} \\
0 &= D_{\mu \alpha p} = 2(d_p - d_3) F^\mu_{\alpha p}
\endaligned
\end{equation}
on $\tilde U$, for all $e,p,\alpha,\mu$;  and from GD(3.59)
\begin{equation}\label{GD359}
\aligned
0 &= E_{p\mu e} = -(e_e + e_3) F^\mu_{p e} \\
0 &= E_{p\mu \alpha} = 2(d_p - d_3) F^\mu_{\alpha p}
\endaligned
\end{equation}
on $\tilde U$, for all $e,p,\alpha,\mu$.
In summary, we have
\begin{equation}\label{U9}
0= (e_e + e_3) F^\alpha_{pe}, \quad 0= (e_e+e_3) F^\mu_{pe}, \quad 0=
(e_e + e_3) F^\mu_{\alpha e} 
\end{equation}
on $\tilde U$ for all $e,p,\alpha, \mu$, and
\begin{equation}\label{U10}
0 = (d_p-d_3) F^\mu_{pe}, \quad 0 = (d_p - d_3) F^\mu_{\alpha p},
\quad 0 = (d_p-d_3) F^\alpha_{pe}
\end{equation}
on $\tilde U$, for all $e,p,\alpha,\mu$.  In the present
proof, equations~\eqref{U10} are not needed, but we record them here
for use in the proof of Theorem~\ref{MT3} below. 
For each $e \in \1$, define
the analytic function on $U$
\begin{equation}\label{U11a}
A_e = \sum_\alpha |v_{e\alpha}|^2 + \sum_\mu |v_{e\mu}|^2 =
2\sum_{p,\alpha} (F^\alpha_{pe})^2 + 2\sum_{p,\mu}(F^\mu_{pe})^2 + 8
\sum_{\alpha, \mu} (F^\mu_{\alpha e})^2
\end{equation}
Let
\begin{equation}\label{U11b}
n(A_e) = \{ x\in U: A_e(x) >0\}
\end{equation}
an open subset of $U$.  Let
$\partial n(A_e)$ be the boundary of $n(A_e)$ in $U$.  If
\begin{equation}\label{U11c}
\tilde U \cap n(A_e) = \emptyset
\end{equation}
for all $e \in \1$, then Theorem~\ref{TS2} applies and we
conclude that $\lambda$ is reducible on $\tilde U$.  The proof of
Theorem~\ref{TS2} shows that 
$e_3$ is constant on $\tilde U$ in this case, by~\eqref{S16a} and the
fact that $\omega^0_0 = 0$ now.
If
\begin{equation}\label{U11d}
\tilde U \cap n(A_c) \neq \emptyset
\end{equation}
for some $c \in \1$, then there exists a point
\begin{equation}\label{U11e}
x \in \left( \tilde U \cap (\cup_e n(A_e))\right) \setminus \cup_e
\partial n(A_e)
\end{equation}
and for such a point there exists a connected open neighborhood $W$ of
$x$ such that $W\subset \tilde U$ and for every $e \in \1$, either
$A_e$ is identically zero on $W$ or $A_e$ is always positive on $W$.
Let
\begin{equation}\label{U11f}
\widetilde \1' = \{a\in \1 : \text{ $A_a = 0$ on $W$}\}
\end{equation}
and let $\1''$ be the complement of $\widetilde \1'$ in $\1$.  Thus,
\begin{equation}\label{U11g}
\1'' = \{ c \in \1 : \text{ $A_c >0$ on $W$}\}
\end{equation}
and $\1'' \neq \emptyset$.
If $c \in \1''$, then for each $x\in W$, there 
exists $p$, $\alpha$, or $\mu$ such that $F^\alpha_{pc}(x) \neq 0$ or
$F^\mu_{pc}(x) \neq 0$ or $F^\mu_{\alpha c}(x) \neq 0$.  Therefore,
\begin{equation}\label{U12}
e_c = -e_3
\end{equation}
on $W$, for any $c \in \1''$, by~\eqref{U9}.  
We have
\begin{equation}\label{U13}
e_3 >0
\end{equation}
on $W$, because if $c \in \1''$, then by~\eqref{3.36A1},
\begin{equation}\label{U14}
0 < A_c = 2m_3(e_3 - e_c) = 4m_3e_3
\end{equation}
on $W$, by~\eqref{U12}.  We next prove that
\begin{equation}\label{U15}
\omega^a_c = 0
\end{equation}
on $W$, for all $a \in \widetilde\1'$ and $c \in \1''$.  To do
this, we 
observe that $E_{ac} = 0$ on $W$, for all $a \in \widetilde
\1'$ and all 
$c \in \1''$, since $E_1$ is diagonalized on $W$.  Therefore,
using~\eqref{U12}, we have
\begin{equation}\label{U16}
\aligned
\sum_j E_{acj} \theta^j &= dE_{ac} +2 E_{ac} \omega^0_0 - \sum_e E_{ec}
\omega^c_a - \sum_e E_{ae} \omega^e_c \\
&= (e_a-e_c) \omega^c_a = 2e_3 \omega^c_a
\endaligned
\end{equation}
on $W$.   But GD(3.50) with~\eqref{U4} and the definition of
$\widetilde\1'$ imply that
\begin{equation}\label{U17}
E_{acj} = 0
\end{equation}
on $W$, for all $a\in \widetilde\1'$, $c \in \1''$, and all $j$.
Then~\eqref{U15} follows from~\eqref{U13},~\eqref{U16},
and~\eqref{U17}.  In addition, $e_3$ must be constant on $W$.
In fact, from GD(3.52) it is seen that $e_{3e} = 0$, $e_{3p} = 0$, and
$e_{3\mu}=0$ on $W$.  If $a \in \widetilde\1'$, then $e_3=e_a$
on $W$, so $e_{3\alpha} = e_{a\alpha} = E_{aa\alpha}=0$ by
GD(3.50).   Hence, $e_{3j} = 0$ on $W$ for all $j$, and so
\begin{equation}\label{U18}
de_3 = de_3 + 2e_3 \omega^0_0 = \sum_j e_{3j} \theta^j =0
\end{equation}
on $W$.

The rest of the proof is now similar to the last part of the proof
of Theorem~\ref{TS2}.  For each $x \in W$,
let $V(x)$ be the subspace of $\R^{n+3}_2$ defined by the span of the
vectors 
\begin{equation}\label{U19}
Y_0,\; Y_c,\; Y_p,\; Y_\alpha,\; Y_\mu,\; Y_{n+2},\; e_3 Y_1 - Y_{n+1}
\end{equation}
at the point $x \in W$, for all $c\in \1''$, $p,\alpha,\mu$ .
Let $V$ be the span of
$V(x)$ for all $x \in W$.  
We want to prove that
$V$ is a subspace of codimension $m+1$, where $m\geq 1$ is the
cardinality of $\widetilde\1'$.  Because the codimension of
$V(x)$ is $m+1$ for any $x \in W$, we will obtain our result if we
prove that $V(x)$ is constant on $W$.
This will be true if we show that the
derivatives of the vectors spanning $V(x)$ are zero modulo $V(x)$, for
every $x \in W$.  This follows because
\begin{equation}\label{U20}
\omega^a_p=0, \quad \omega^a_\alpha = 0, \quad \omega^a_\mu =0
\end{equation}
on $W$, for all $a \in \widetilde\1'$, by~\eqref{GD(3.13)} and the
definition of $\widetilde\1'$. 
\end{proof}

\begin{corollary}\label{Cor22}
Let $\lambda:M \to \Lambda$ be an irreducible proper
Dupin 
hypersurface for which $r=-1$, $m_1=m_2$, and $m_3=m_4$.  If $Y:U\to
G$ is a second order frame field such that $d\omega^0_0=0$ 
on $U$, then
\begin{equation}\label{U5cor}
\1' = \emptyset
\end{equation}
\end{corollary}

\section{One pair of multiplicities is 1}\label{Sec8}

Assume now that $\lambda:M^{n-1} \to \Lambda^{2n-1}$ is the Legendre
lift of a Dupin hypersurface with four principal curvature, constant
Lie curvature $r=-1$ and multiplicities $m_1=m_2 \geq 2$ and
$m_3=m_4=1$.  For these multiplicities, the index sets $\3$ and $\4$
consist of one element each
\[
\3 = \{2m_1+2\},\; \4 = \{2m_1+3\}
\]
It is convenient to continue writing $\alpha$ and $\mu$ for these
values, respectively.  In addition, $D_3$, $D_4$, $E_3$, and $E_4$ are
automatically scalar matrices in this case, for any second order frame
field.

\begin{proposition}\label{PM1}  If $m_1=m_2 \geq 2$, $m_3=m_4=1$, and
$r=-1$, then for any point in $M$
there exists a second order frame
field $Y:U \to 
G$ about the point, for which
\begin{equation}\label{4}
d_1 = e_2
\end{equation}
For any such frame field,
\begin{equation}\label{5}
F^\alpha_{pa} F^\mu_{pa} F^\mu_{\alpha a} = 0 = F^\alpha_{pa}
F^\mu_{pa} F^\mu_{\alpha p}
\end{equation}
on $U$, for all $ p,a$; and
\begin{equation}\label{M0}
d_{1a} = 0 = R_a, \quad d_{1p} = 0 = R_p
\end{equation}
on $U$, for all $a,p$.
\end{proposition}

\begin{proof} 
Let $Y:U\to G$ be a second order frame field.  Then $D_1 = d_1I_{m_1}$
and $E_2= e_2I_{m_1}$ are scalar matrices, by the first two equations
in GD(3.42).  A second order frame change~\eqref{GD(3.14)} of the form
$\tilde Y = Y a(I_2,I,0,sL)$ has
\begin{equation}\label{M1}
\tilde e_2 - \tilde d_1 =e_2 - d_1 + 2s
\end{equation}
by GD(3.32).  Taking $s = (d_1-e_2)/2$, we obtain a second order frame
field for which 
\begin{equation}\label{M2}
d_1=e_2
\end{equation}
We assume this done for our frame $Y$.  Setting $a=b=c$ in GD(3.62i),
we find that
\begin{equation}\label{M3}
\sum_p D_{p\alpha} F^\alpha_{pc} + \sum_p E_{p\mu} F^\mu_{pc} = 4
\sum_p F^\alpha_{pc} F^\mu_{pc} F^\mu_{\alpha c}
\end{equation}
By GD(3.51i), for each $p=q$ we have
\begin{equation}\label{M4}
{e_2}_c = E_{ppc} = R_c + 2E_{p\mu} F^\mu_{pc} + 2 D_{p\alpha}
F^\alpha_{pc} 
\end{equation}
which shows that
\begin{equation}\label{M5}
E_{p\mu} F^\mu_{pc} +  D_{p\alpha}
F^\alpha_{pc}
\end{equation}
is independent of $p \in \2$.  Therefore,~\eqref{M3} becomes
\begin{equation}\label{M6}
D_{p\alpha} F^\alpha_{pc} + E_{p\mu} F^\mu_{pc} = \frac4{m_1}
\sum_q F^\alpha_{qc} F^\mu_{qc} F^\mu_{\alpha c}
\end{equation}
for all $p \in \2$.
Using GD(3.46i) with $a=b=c$, and using~\eqref{M3}, we find
\begin{equation}\label{M7}
{d_1}_c=D_{ccc} = -R_c
\end{equation}
By~\eqref{M2}, ${d_1}_c = {e_2}_c$, so
combining~\eqref{M4} and~\eqref{M7} and substituting the result
into~\eqref{M3} gives
\begin{equation}\label{M8}
\-R_c = \frac 4{m_1} \sum_p F^\alpha_{pc}F^\mu_{pc}F^\mu_{\alpha c}
\end{equation}
Multiply $F^\alpha_{pcc}$ by $2F^\alpha_{pc}$ and $F^\mu_{pcc}$
by $2F^\mu_{pc}$ in GD(3.37) and GD(3.38), respectively, subtract the
latter from the former and use GD(3.36i), to get
\begin{equation}\label{M9}
\aligned
-d_{1c} &= ((F^\alpha_{pc})^2 - (F^\mu_{pc})^2)_c = 2F^\alpha_{pc}
F^\alpha_{pcc} -2 F^\mu_{pc}F^\mu_{pcc} \\
&= 2(F^\alpha_{pc} D_{p\alpha} + F^\mu_{pc} E_{p\mu}) -12
F^\mu_{pc}F^\mu_{\alpha c} F^\alpha_{pc}
\endaligned
\end{equation}
Therefore, for each $c \in \1$,
\begin{equation}\label{M10}
F^\mu_{pc}F^\mu_{\alpha c} F^\alpha_{pc}
\end{equation}
is independent of $p$, by~\eqref{M5} and the fact that $d_{1c}$ is
independent of $p$.  We show now that this implies~\eqref{5}.  At any
point of $U$, let $\mathcal V$ denote the vector subspace spanned by
the $Y_p$.  This subspace is invariant under a change of
frame~\eqref{GD(3.14)}.  On $\mathcal V$, for fixed
$c$,  define the bilinear form
\begin{equation}\label{M11}
\aligned
S:\mathcal V\times \mathcal V &\to \R \\
S(\sum_p u^p Y_p, \sum_q v^q Y_q) &= \sum_{p,q} F^\mu_{pc}F^\mu_{\alpha
c} F^\alpha_{qc} u^p v^q
\endaligned
\end{equation}
By GD(3.10), $S$ does not depend on the choice of $Y$.
Then~\eqref{M10} is the value of $S(Y_p,Y_p)$, so this value is
independent of $p$.  Let
\begin{equation}\label{M11a}
K = S(Y_p,Y_p)
\end{equation}
for every $p \in \2$, an analytic function on $U$.  If $t$ and $s$ are
nonzero real numbers such 
that $t^2 + s^2 =1$, and if $p \neq q$, then replacing $Y_p$ and $Y_q$
with $tY_p + sY_q$ and $-sY_p + tY_q$ gives another allowable frame
$Y$, so $S$ has the same value on each:
\begin{equation}\label{M11b}
S(tY_p + sY_q,tY_p + sY_q) = S(-sY_p + tY_q,-sY_p + tY_q)
\end{equation}
so
\begin{equation}\label{M11c}
K + 2ts S(Y_p,Y_q) = K - 2st S(Y_p,Y_q)
\end{equation}
Therefore,
\begin{equation}\label{M11d}
S(Y_p,Y_q) = 0
\end{equation}
whenever $p \neq q$, so
$S$ is a multiple of the inner product on $\mathcal V$;
that is,
\begin{equation}\label{M12}
S(\sum_p u^pY_p, \sum_p u^p Y_p) = K \sum_p (u^p)^2
\end{equation}
for some function $K$ on $U$.  On the other hand,
\begin{equation}\label{M13}
S(\sum_p u^pY_p, \sum_p u^p Y_p) = F^\mu_{\alpha c}( \sum_p
F^\mu_{pc}u^p) (\sum_p F^\alpha_{pc} u^p)
\end{equation}
is the product of two linear polynomials in the $u^p$.  Such a
factorization of a sum of two or more squares is impossible over the reals unless
$K$ is identically zero on $U$.  Therefore, the first
equation~\eqref{5} must hold at 
every point of $U$.
The proof of the second equation~\eqref{5} is done in the same way with
the roles of $a$ and $p$ reversed.
\end{proof}


\begin{proposition}\label{PM2} If $m_1 = m_2 \geq 2$, $m_3=m_4=1$, and
$r=-1$, and if
$Y:U\to G$ is a second order frame
field for which $d_1+e_2=0$ on $U$, then
$d\omega^0_0 = 0$ on $U$. 
\end{proposition}

\begin{proof} 
As seen in the proof of Proposition~\ref{PM1}, we may assume $Y$
chosen so that $d_1=e_2$ on $U$, in which case the hypothesis implies
that $d_1=e_2=0$ on $U$.  Then GD(3.36i) says
\begin{equation}\label{M15}
(F^\alpha_{pa})^2 = (F^\mu_{pa})^2
\end{equation}
on $U$, for all $a,p$.  Then adding together GD(3.36ii) and
GD(3.36iv), respectively, subtracting GD(3.36iii) from GD(3.36v), we find
\begin{equation}\label{M16}
e_3 = -e_4, \quad d_3 = d_4
\end{equation}
respectively, on $U$.  Substituting~\eqref{M15} into~\eqref{5} we find
\begin{equation}\label{M16a}
F^\mu_{pa} F^\mu_{\alpha a} = 0
\end{equation}
on $U$, for all $a,p$, and
\begin{equation}\label{M16b}
F^\mu_{pa}F^\mu_{\alpha p} = 0
\end{equation}
on $U$, for all $a,p$.  These equations are true even for $a\neq b$
and $p \neq q$, as follows.  In fact, take the covariant derivative
of~\eqref{M15} to get
\begin{equation}\label{M16c}
F^\alpha_{pa}F^\alpha_{paj} = F^\mu_{pa}F^\mu_{paj}
\end{equation}
for all $a,p,j$.  Take $j=b \neq a$ and use GD(3.37) and GD(3.38) to
get
\begin{equation}\label{M16d}
F^\alpha_{pa}(-F^\mu_{pa}F^\mu_{\alpha b} - 2F^\mu_{pb}F^\mu_{\alpha
a}) = F^\mu_{pa}(F^\alpha_{pa} F^\mu_{\alpha b} + 2 F^\alpha_{pb}
F^\mu_{\alpha a})
\end{equation}
on $U$.  This, together with~\eqref{5}, yields
\begin{equation}\label{M16e}
F^\alpha_{pa} F^\mu_{pa}F^\mu_{\alpha b} + F^\alpha_{pa} F^\mu_{pb}
F^\mu_{\alpha a} + F^\alpha_{pb} F^\mu_{pa}F^\mu_{\alpha a} =0
\end{equation}
on $U$, for all $a,b,p$.
By~\eqref{M16a} and~\eqref{M15}, the middle term is 0 on $U$,
while~\eqref{M16a} implies that the third term is 0 on $U$.
Thus~\eqref{M16e} is equivalent to
$F^\alpha_{pa}F^\mu_{pa}F^\mu_{\alpha b}=0$ on $U$, for all $a,b,p$, which
by~\eqref{M15}, is equivalent to
\begin{equation}\label{M17}
F^\mu_{pa} F^\mu_{\alpha b} = 0
\end{equation}
on $U$, for all $a, b, p$.  
In the same way, if we take $j=q\neq p$ in~\eqref{M16d}, and use
GD(3.37) and GD(3.38), we get
\begin{equation}\label{M17a}
F^\alpha_{pa}(F^\mu_{pa}F^\mu_{\alpha q} + 2 F^\mu_{\alpha p}
F^\mu_{qa}) = F^\mu_{pa}(-F^\alpha_{pa} F^\mu_{\alpha q} - 2
F^\mu_{\alpha p} F^\alpha_{qa})
\end{equation}
on $U$, which simplifies to
\begin{equation}\label{M17b}
F^\alpha_{pa}F^\mu_{pa}F^\mu_{\alpha q} + F^\alpha_{pa} F^\mu_{qa}
F^\mu_{\alpha p} + F^\alpha_{qa}F^\mu_{pa} F^\mu_{\alpha p} = 0
\end{equation}
on $U$, for all $a, p\neq q$.  By~\eqref{M16b} and~\eqref{M15}, the
middle term is 0 on $U$, 
while~\eqref{M16b} implies that the third term is 0 on $U$.
Thus~\eqref{M17b} is equivalent to
$F^\alpha_{pa}F^\mu_{pa}F^\mu_{\alpha q}=0$ on $U$, for all $a,p \neq q$,
which by~\eqref{M15}, is equivalent to
\begin{equation}\label{M18}
F^\mu_{pa}F^\mu_{\alpha q} = 0
\end{equation}
on $U$, for all $a$, $p\neq q$.  By~\eqref{M16b}, this holds also for
$p=q$. 

For a function $f:U \to \R$, let 
\begin{equation}\label{M19}
n(f) = \{ x\in U: f(x) \neq 0 \}
\end{equation}
Let $\spt(f)$ denote the closure of $n(f)$ in $U$, namely, the support
of $f$ in $U$.  Then~\eqref{M15} implies that
\begin{equation}\label{M20}
n(F^\alpha_{pa}) = n(F^\mu_{pa})
\end{equation}
for all $a,p$.  Define open subsets $U_2$ and $U_3$ of $U$ by
\begin{equation}\label{M20a}
U_2 = \cup_{a,p} n(F^\mu_{pa}), \quad U_3 = \cup_a n(F^\mu_{\alpha a})
\end{equation}
Then~\eqref{M17} implies that
\begin{equation}\label{M21}
U_2 \cap U_3 = \emptyset 
\end{equation}
and \eqref{M18} implies that
\begin{equation}\label{M22}
U_2 \cap (\cup_q n(F^\mu_{\alpha q})) = \emptyset 
\end{equation}
Let $U_1 \subset U$ be the open subset of $U$ defined by 
\begin{equation}\label{M23}
U_1 = U \setminus ((\cup_{a,p} \spt(F^\mu_{pa})) \cup (\cup_{a}
\spt(F^\mu_{\alpha a})))
\end{equation}
The closure of $U_2 \cup U_3$ in $U$ is clearly the complement of
$U_1$ in $U$. 

On $U_1$, the functions $f$, $g$, and $h$ in~\eqref{3} are identically
zero, by~\eqref{M20} and the definition of $U_1$.  Hence,
by~\eqref{S15}  in the proof of Theorem~\ref{TS2}, we have
$d\omega^0_0 = 0$ on $U_1$.

On $U_2$ we have 
\begin{equation}\label{M23a}
F^\mu_{\alpha a} = 0, \quad F^\mu_{\alpha p} = 0
\end{equation}
for all $a,p$, by~\eqref{M21} and~\eqref{M22}, respectively.
Then by~\eqref{GD(3.34)}, $F^\mu_{\alpha a j} = 0$ on $U_2$, for all
$j$.  By GD(3.39), we have
\begin{equation}\label{M24}
0 = F^\mu_{\alpha a a} = -\frac 12 D_{\mu \alpha} - \frac 12 \sum_p
F^\mu_{pa} F^\alpha_{pa} + \frac 12 \sum_p F^\mu_{pa} F^\alpha_{pa} =
-\frac12 D_{\mu \alpha}
\end{equation}
on $U_2$.  By GD(3.39) and~\eqref{M23a} we have
\begin{equation}\label{M26}
0 = F^\mu_{\alpha a \alpha} = D_{\mu a} + 3 \sum_p F^\mu_{\alpha p}
F^\alpha_{pa} = D_{\mu a}
\end{equation}
on $U_2$, for all $a$.  By GD(3.39) and~\eqref{M23a}, we have
\begin{equation}\label{M27}
0 = F^\mu_{\alpha a \mu} = D_{a\alpha} + 3 \sum_p F^\mu_{\alpha p}
F^\mu_{pa} = D_{a \alpha}
\end{equation}
on $U_2$, for all $a$.  Similarly, by the second equation
in~\eqref{M23a}, we have $F^\mu_{\alpha p j} = 0$ on $U_2$, for all
$p,j$.  From GD(3.40) and~\eqref{M23a} we get
\begin{equation}\label{M28}
0 = F^\mu_{\alpha p \alpha} = - E_{p\mu} - 3 \sum_a F^\alpha_{pa}
F^\mu_{\alpha a} = -E_{p\mu}
\end{equation}
on $U_2$, for all $p$.  Similarly,
\begin{equation}\label{M29}
0 = F^\mu_{\alpha p \mu} = -D_{p\alpha} - 3 \sum_a F^\mu_{\alpha a}
F^\mu_{pa} = -D_{p\alpha}
\end{equation}
on $U_2$, for all $p$.  Finally, we want to show that $D_{pa} =0$ on
$U_2$, for all $a,p$.  By GD(3.37)  and the fact that
$F^\mu_{\alpha a}=0$ and $F^\mu_{\alpha p}=0$ at every point of $U_2$,
for all $a,p$, we have
\begin{equation}\label{M30}
F^\alpha_{pa\alpha} = -D_{pa}
\end{equation}
on $U_2$, for all $a,p$.  In the same way using GD(3.38), we have
\begin{equation}\label{M31}
F^\mu_{pa\mu} = -D_{pa}
\end{equation}
on $U_2$, for all $a,p$.  Taking $j=\alpha$ in~\eqref{M16c}, we find
\begin{equation}\label{M32}
F^\alpha_{pa} F^\alpha_{pa\alpha} = F^\mu_{pa} F^\mu_{pa\alpha}
\end{equation}
on $U_2$, for all $a,p$. Substitute~\eqref{M30} into this to get
\begin{equation}\label{M33}
-F^\alpha_{pa}D_{pa} = F^\mu_{pa} F^\mu_{pa\alpha}
\end{equation}
on $U_2$, for all $a,p$.  Similarly, taking $j=\mu$ in~\eqref{M16c},
we get
\begin{equation}\label{M34}
F^\alpha_{pa} F^\alpha_{pa\mu} = F^\mu_{pa} F^\mu_{pa\mu}
\end{equation}
on $U_2$, for all $a,p$.  Substitute~\eqref{M31} into this to get
\begin{equation}\label{M35}
F^\alpha_{pa}F^\alpha_{pa\mu} = - F^\mu_{pa} D_{pa}
\end{equation}
on $U_2$, for all $a,p$.  By GD(3.41) we have $F^\alpha_{pa\mu} = -
F^\mu_{pa\alpha}$, so substitute this into~\eqref{M35} to get
\begin{equation}\label{M36}
F^\alpha_{pa} F^\mu_{pa\alpha} = F^\mu_{pa} D_{pa}
\end{equation}
on $U_2$, for all $a,p$.  Now multiply~\eqref{M33} by $F^\mu_{pa}$ to
get
\begin{equation}\label{M37}
-F^\alpha_{pa}F^\mu_{pa}D_{pa} = (F^\mu_{pa})^2 F^\mu_{pa\alpha}
\end{equation}
on $U_2$, for all $a,p$.  Multiply~\eqref{M36} by $F^\alpha_{pa}$ to
get 
\begin{equation}\label{M38}
F^\alpha_{pa} F^\mu_{pa} D_{pa} = (F^\alpha_{pa})^2 F^\mu_{pa\alpha}
\end{equation}
on $U_2$, for all $a,p$.
By~\eqref{M15}, we see that the right sides of~\eqref{M37}
and~\eqref{M38} are the same.  Therefore the left sides must be
equal at every point of $U_2$, but they are negatives of each other,
so we get, again using~\eqref{M15}
\begin{equation}\label{M39}
F^\mu_{pa} D_{pa} =0, \quad F^\alpha_{pa} D_{pa} =0
\end{equation}
on $U_2$, for all $a,p$.

The next step is to show 
\begin{equation}\label{M40}
\omega^0_{n+1} = 0
\end{equation}
on $U_2$.  By~\eqref{GD(3.45)}, this is equivalent to showing $R_j=0$
on $U_2$, for all $j$.  Now $d_1=0$ on $U_2$ implies that all of its
covariant derivatives $d_{1j} = 0$ on $U_2$, since $d d_1 + 2d_1
\omega^0_0 = \sum_j d_{1j} \theta^j$ defines its covariant
derivatives.   Using the fact that $D_{\mu\alpha}$, $D_{\mu a}$,
$D_{a\alpha}$, $E_{p\mu}$, $D_{p\alpha}$, and all their covariant
derivatives are zero on $U_2$, we  
find the consequences of $d_{1j}=0$ on $U_2$ to be as follows.
By GD(3.46i) and~\eqref{M23a},
\begin{equation}\label{M41}
0 = d_{1a} = D_{bba} = -R_a
\end{equation}
on $U_2$, for all $a$.  By GD(3.46ii),
\begin{equation}\label{M42}
0 = d_{1p} = D_{aap} = -R_p
\end{equation}
on $U_2$, for all $p$.  By GD(3.46iii) and~\eqref{M39},
\begin{equation}\label{M43}
0 = d_{1\alpha} = D_{aa\alpha} = -R_\alpha
\end{equation}
on $U_2$.  By GD(3.46iv) and~\eqref{M39},
\begin{equation}\label{M44}
0 = d_{1\mu} = D_{aa\mu} = -R_\mu 
\end{equation}
on $U_2$.  This completes the proof of~\eqref{M40}.

The next step is to show that the covariant derivatives $e_{3j}$ of
$e_3$ satisfy
\begin{equation}\label{M45}
e_{3j} = 0
\end{equation}
on $U_2$, where by~\eqref{GD(3.43)} these are defined by
\begin{equation}\label{M46}
d e_3 + 2 e_3 \omega^0_0 = \sum_j e_{3j} \theta^j
\end{equation}
To prove~\eqref{M45}, we use the equations GD(3.52) to get
\begin{equation}\label{M47}
e_{3a} = E_{\alpha\alpha a} = R_a =0
\end{equation}
on $U_2$, for all $a$,
\begin{equation}\label{M48}
e_{3p} = E_{\alpha\alpha p} = R_p = 0
\end{equation}
on $U_2$, for all $p$,
\begin{equation}\label{M49}
e_{3\alpha} = E_{\alpha\alpha\alpha} = R_\alpha - \frac 6{m_1}
\sum_{a,p} D_{pa}F^\alpha_{pa} + \frac{24}{m_1} \sum_a F^\mu_{\alpha
a} F^\alpha_{pa} F^\mu_{\alpha p} = 0
\end{equation}
on $U_2$, by~\eqref{M39} and~\eqref{M23a}, and
\begin{equation}\label{M49a}
e_{3\mu} = E_{\alpha \alpha \mu} = R_\mu + D_{\mu\alpha\alpha}=0
\end{equation}
on $U_2$, since~\eqref{M24} implies that $D_{\mu\alpha\alpha}=0$ on
$U_2$. This completes the proof
of~\eqref{M45}, and therefore~\eqref{M46} becomes
\begin{equation}\label{M50}
d e_3 + 2e_3 \omega^0_0 = 0
\end{equation}
on $U_2$.  This implies that on the subset of $U_2$ where $e_3 \neq
0$, the 1-form
\begin{equation}\label{M51}
\omega^0_0 = -\frac12 \frac 1{e_3} de_3
\end{equation}
and therefore $\omega^0_0$ is closed on this subset of $U_2$.  Thus,
if $e_3$ is never zero on $U_2$, then
$\omega^0_0$ is closed on all of $U_2$.  We now prove that
\begin{equation}\label{M52}
e_3(x) \neq 0
\end{equation}
for every $x \in U_2$.  Let $x \in U_2$.  We may assume that our
second order frame $Y$ on $U$ diagonalizes $E_1$ at $x$.  It is easily
checked that if the vectors $Y_a$ are changed by an orthogonal matrix,
then the set $U_2$ does not change, and $e_3$ is unchanged.  Let
\begin{equation}\label{M53}
\1'(x) = \{ a\in \1 : F^\alpha_{pa}(x) = 0, \;\; \forall p \}
\end{equation}
and let
\begin{equation}\label{M54}
\1''(x) = \1 \setminus \1'(x)
\end{equation}
We want to show that if the set $\1$ is arranged with the indices in
$\1'(x)$ first followed by the indices in $\1''(x)$, then
\begin{equation}\label{M55}
E_1(x) = e_3(x) \begin{pmatrix} I_m & 0 \\ 0 & -I_{m_1 - m}
\end{pmatrix}
\end{equation}
where $m$ is the cardinality of $\1'(x)$.  To prove this, suppose that
$a \in \1'(x)$.  Then $F^\alpha_{pa}(x) = 0$ for all $p$, and of
course $F^\mu_{\alpha a}=0$ on all of $U_2$, so GD(3.36ii) implies
that (since $d_1=0$)
\begin{equation}\label{M56}
0 = -e_a(x) + e_3(x)
\end{equation}
as claimed.  Next suppose that $a \in \1''(x)$, so that $F^\mu_{pa}(x)
= \pm F^\alpha_{pa}(x) \neq 0$, for some $p$.  Using GD(3.59i) and the
fact that $E_{p\mu}=0$ on $U_2$, we have
\begin{equation}\label{M57}
0=E_{p\mu a}(x) = (e_4(x) - e_a(x))F^\mu_{pa}(x)
\end{equation}
which implies that $e_a(x) = e_4(x) = -e_3(x)$, by~\eqref{M16}.  This
completes the proof of~\eqref{M55} and allows us to prove~\eqref{M52}
as follows.  For our given $x \in U_2$, GD(3.36ii) implies that
\begin{equation}\label{M58}
\sum_p (F^\alpha_{pa}(x))^2 = -e_a(x) + e_3(x)
\end{equation}
for all $a$.
But~\eqref{M55} implies that $e_a(x) = \pm e_3(x)$, so if $e_3(x)=0$,
then $e_a(x) =0$ as well, and we must conclude from~\eqref{M58} that
$F^\alpha_{pa}(x) =0$, for all $a,p$, contradicting the fact that $x
\in U_2$.  We conclude that $e_3$ is never zero on $U_2$, and
therefore $\omega^0_0$ is closed on $U_2$, so in particular, 
\begin{equation}\label{M59}
D_{pa} = 0
\end{equation}
on $U_2$, for all $a,p$.

Finally, it remains to prove that $d\omega^0_0=0$ on $U_3$.
To do this, we first observe that $F^\alpha_{pa} = 0 = F^\mu_{pa}$ on
$U_3$, for every $a,p$.  It follows that their covariant derivatives
$F^\alpha_{paj}=0= F^\mu_{paj}$ on $U_3$, for all $a,p,j$.  By
GD(3.37i), then
\begin{equation}\label{60}
0 = F^\alpha_{paa} = D_{p\alpha}
\end{equation}
on $U_3$, for all $a,p$.  By GD(3.37ii),
\begin{equation}\label{M61}
0 = F^\alpha_{pap} = -D_{a\alpha}
\end{equation}
on $U_3$, for all $a,p$.  By GD(3.37iii)
\begin{equation}\label{M62}
0 = F^\alpha_{pa\alpha} = -D_{pa} + 2F^\mu_{\alpha p} F^\mu_{\alpha a}
-2 F^\mu_{\alpha p} F^\mu_{\alpha a} = -D_{pa}
\end{equation}
on $U_3$, for all $a,p$.  By GD(3.38i)
\begin{equation}\label{M63}
0 = F^\mu_{paa} = -E_{p\mu}
\end{equation}
on $U_3$, for all $a,p$.  By GD(3.38ii)
\begin{equation}\label{M64}
0 = F^\mu_{pap} = D_{\mu a}
\end{equation}
on $U_3$, for all $a,p$.  It remains to prove that $D_{\mu\alpha}=0$
on $U_3$.  To do this, we first observe that there is a dense open
subset $V$
of $U_3$ on which we may make an analytic change of frame field 
\begin{equation}\label{M64a}
\tilde Y = Y a(I_2, B,0,0)
\end{equation}
where $B = B_1 \oplus I_{m_1} \oplus 1 \oplus 1$, with $B_1:V \to
O(m_1)$ an analytic map, such that
$\tilde E_1$ is diagonal at each point (see, for
example,~\cite{Singley}). For such a change of frame, $\tilde
\omega^0_0 = \omega^0_0$ and the sets $\tilde U_2 = U_2$ and $\tilde
U_3 = U_3$.  So, without loss of generality, we may assume $\tilde Y =
Y$ on $V$ and all of the above properties of $Y$ on $U_3$ remain true
on $V$.
For a point $x \in V$, let
\begin{equation}\label{M65}
\1'(x) = \{a \in \1: F^\mu_{\alpha a}(x) = 0 \}
\end{equation}
and let its complement be
\begin{equation}\label{M66}
\1''(x) = \1 \setminus \1'(x)
\end{equation}
If $a \in \1'(x)$, then GD(3.36ii) implies that
\begin{equation}\label{M67}
0 = -e_a(x) + e_3(x)
\end{equation}
If $a \in \1''(x)$, then by~\eqref{M61} and GD(3.54iv) and
using~\eqref{M16}, we have
\begin{equation}\label{M68}
0 = D_{a\alpha \mu} = (-2e_a-e_3-d_3+d_4+e_4)F^\mu_{\alpha a} =
-2(e_a+e_3)F^\mu_{\alpha a}
\end{equation}
at $x$, and so $F^\mu_{\alpha a}(x) \neq 0$ implies that 
\begin{equation}\label{M69}
e_a(x) = -e_3(x)
\end{equation}
We next show that $e_3(x) \neq 0$.  Since $x \in V \subset U_3$, the
set $\1''(x)$ is non-empty.  For $a \in \1''(x)$, GD(3.36ii)
and~\eqref{M69} give
\begin{equation*}
4(F^\mu_{\a a})^2 = -e_a(x) + e_3(x) = 2e_3(x)
\end{equation*}
The left side of this equation is non-zero for $a \in \1''(x)$, and
thus $e_3(x) \neq 0$.
Therefore, at every point of $V$, with this frame $ Y$ the
matrix 
$E_1$ has the form~\eqref{M55}, where $m = m(x)$ is the cardinality of
$\1'(x)$.  Since $e_3(x) \neq 0$, this shows that on
each connected component of $V$, $\1'$ and $\1''$ must be constant and
equations~\eqref{M67} and~\eqref{M69} hold at every point of this component.  Therefore, we may use
GD(3.50i) to find that the covariant derivative of $e_a$, for any $a \in
\1''$, is
\begin{equation}\label{M70}
e_{aa} = E_{aaa} = 6D_{\mu\alpha}F^\mu_{\alpha a}
\end{equation}
on $V$, while by~\eqref{M69} and GD(3.52i)
\begin{equation}\label{M71}
-e_{aa} = e_{3a} = E_{\alpha\alpha a} = R_a +
2D_{\mu\alpha}F^\mu_{\alpha a} = 
2D_{\mu\alpha}F^\mu_{\alpha a} 
\end{equation}
on $V$, by~\eqref{5} and~\eqref{M8}.  Adding~\eqref{M71}
and~\eqref{M70}, we get
\begin{equation}\label{M72}
0 = 8D_{\mu\alpha}F^\mu_{\alpha a}
\end{equation}
on this connected component, for any $a \in \1''$.  Therefore,
\begin{equation}\label{M73}
D_{\mu\alpha}=0
\end{equation}
on each connected component of $V$, therefore on all
of $V$.  It follows from~\eqref{GD(3.33)} that $d\omega^0_0=0$ on $V$,
and therefore
$d\omega^0_0 = 0$ on the closure of $V$ in
$U_3$, which is all of $U_3$.  Therefore, this is true for our
originally chosen frame field $Y$ on $U$, and $d\omega^0_0 = 0$ on
$U_1$, $U_2$, and $U_3$, so it is also zero on the closure of the
union of these three sets in $U$, which is all of $U$.
\end{proof}

\begin{theorem}\label{PM3} Suppose $m_1=m_2 \geq 2$, $m_3=m_4 =1$, and
$r=-1$. 
For any point
of $M$, there exists a
second order frame field $Y:U \to G$ about the point, for which
$d\omega^0_0=0$ on $U$.
\end{theorem}

\begin{proof} 
Given a point of $M$, we know by Proposition~\ref{PM1} that there
exists a second order frame field $Y:U\to G$ about the point for which
$d_1=e_2$ on $U$ and~\eqref{5} holds on $U$, for all $a,p$.  
We will show that $d\omega^0_0 =0$ on $U$, for this frame.  To do
this, we decompose $U$ into a disjoint union of subsets $U_1$ and
$U_2$, where 
\begin{equation}\label{M73a}
U_1 = \{ x\in U: d_1(x) = 0 \}
\end{equation}
a closed subset of $U$, and
\begin{equation}\label{M73b}
U_2 = \{ x \in U: d_1(x) \neq 0\}
\end{equation}
an open subset of $U$.  Let $\stackrel{\circ}{U_1}$ denote the
interior of $U_1$.  Then $d_1=0$ on $\stackrel{\circ}{U_1}$, implies
$d\omega^0_0 = 0$ on $\stackrel{\circ}{U_1}$, by Proposition~\ref{PM2}. 

Next consider the open set $U_2$, where $d_1$ is never zero.
On each connected component of $U_2$,  $d_1$ is either always
positive or always negative.
To be specific, suppose that 
\begin{equation}\label{M74a}
e_2=d_1<0
\end{equation}
on the connected component $\tilde U$ of $U_2$.  Then GD(3.36i) becomes
\begin{equation}\label{M74}
2(F^\alpha_{pa})^2 - 2(F^\mu_{pa})^2 = -(d_1+e_2) = -2d_1 >0
\end{equation}
on $\tilde U$ for all $a,p$, which implies that
\begin{equation}\label{M75}
F^\alpha_{pa}\neq 0
\end{equation}
at every point of $\tilde U$, for all $a,p$.  This combined
with~\eqref{5} 
implies that 
\begin{equation}\label{M76}
F^\mu_{pa}F^\mu_{\alpha a} = 0, \quad F^\mu_{pa} F^\mu_{\alpha p} = 0
\end{equation}
on $\tilde U$, for all $a,p$.  This holds for an arbitrary orthogonal
change 
of the frame vectors $Y_a$ or of the frame vectors $Y_p$.  By
polarization, it follows that
\begin{equation}\label{M77}
F^\mu_{pa}F^\mu_{\alpha b} + F^\mu_{pb} F^\mu_{\alpha a} = 0, \quad
F^\mu_{pa} F^\mu_{\alpha q} + F^\mu_{qa} F^\mu_{\alpha p} = 0
\end{equation}
on $\tilde U$, for all $a,b,p,q$.  Polarization means the following.
For any 
$p$ consider the bilinear form
\begin{equation}\label{M77a}
\Phi_p = \sum_{a,b} F^\mu_{pa}F^\mu_{\alpha b} \theta^a \otimes
\theta^b
\end{equation}
Then~\eqref{M76} says that $\Phi_p(Y_a,Y_a)=0$ for any $a$, and this
is true for any orthogonal change of the frame vectors $Y_a$, so it
must be true that $\Phi_p(v,v)=0$ for any unit vector $v$.  Hence,
$\Phi_p$ is an alternating form and must satisfy
\begin{equation}\label{M77b}
\Phi_p(u,v) + \Phi_p(v,u) = 0
\end{equation}
for all vectors $u,v$.  This implies the first equation
in~\eqref{M77}.  The second equation follows in the same way.

The next step is to prove that
\begin{equation}\label{M78}
F^\mu_{pa} = 0
\end{equation}
on $\tilde U$, for all $a,p$.  We prove this by assuming the contrary
and 
then deducing a contradiction to our assumption that $d_1 \neq 0 $ on
$\tilde U$.  Suppose,
then, there exists $a$ and $p$ and $x \in \tilde U$
such that 
\begin{equation}\label{M78a}
F^\mu_{pa}(x) \neq 0
\end{equation}
Then there is an open set $O$ of $\tilde U$, containing $x$, on
which $F^\mu_{pa}$ is never zero, for this $a,p$.
Hence $F^\mu_{\alpha a}=0$ on $O$, for
this $a$, by~\eqref{M76}, and then 
\begin{equation}\label{M79}
F^\mu_{\alpha b} = 0
\end{equation}
on $O$, for all $b$, by the first equation in~\eqref{M77}.
In the same way we see that
\begin{equation}\label{M80}
F^\mu_{\alpha q} = 0
\end{equation}
on $O$, for all $q$.  These two equations are the same
as~\eqref{M23a}, and in the same way as we derived~\eqref{M24}
through~\eqref{M29}, we get
\begin{equation}\label{M81}
D_{\mu\alpha}= D_{\mu b} = D_{b\alpha} = E_{q\mu} = D_{q\alpha}=0
\end{equation}
on $O$, for all $b,q$.  There exists an open dense subset $\tilde
O$ of $O$ on which we may assume $D_2$ and $E_1$ have been
diagonalized at each point.  Then
\begin{equation}\label{M82}
\aligned
0 &= D_{q\alpha b} = (e_3 - e_2 + e_b) F^\alpha_{qb} \\
0 &= E_{p\mu a} = (e_4 - e_2 - e_a)F^\mu_{pa}
\endaligned
\end{equation}
on $\tilde O$, by GD(3.56) and GD(3.59).  
By~\eqref{M75} we can conclude that
\begin{equation}\label{M83}
e_b = e_2 - e_3
\end{equation}
on $\tilde O$, for all $b$, so by continuity, $E_1 =
(e_2-e_3)I_{m_1}$ is scalar on all of $O$.  But,
by~\eqref{M78a} and the second equation of~\eqref{M82}, we have
\begin{equation}\label{M84}
e_a = e_4 - e_2 = -e_3 -e_2
\end{equation}
on $O$, for the particular value of $a$, therefore for all
$a$, since $E_1$ is scalar on $O$.  Subtracting~\eqref{M84}
from~\eqref{M83}, we conclude that $e_2 
= 0$ on $O$, which contradicts~\eqref{M74a}.  

Therefore,~\eqref{M78} holds on $\tilde U$, for all $a,p$.  Then all
covariant 
derivatives $F^\mu_{paj}=0$ on $\tilde U$, for all $a,p,j$, so GD(3.38)
implies that
\begin{equation}\label{M85}
0 = F^\mu_{pa\mu} = -D_{pa}
\end{equation}
on $\tilde U$, for all $a,p$, and
\begin{equation}\label{M86}
\aligned
0 &= F^\mu_{pap} = D_{\mu a} - 3F^\alpha_{pa} F^\mu_{\alpha p} \\
0 &= F^\mu_{paa} = -E_{p\mu} + 3F^\alpha_{pa}F^\mu_{\alpha a}
\endaligned
\end{equation}
on $\tilde U$, for all $a,p$.  These two formulas hold for any
orthogonal change of the vectors $Y_a$ or $Y_p$ on $\tilde U$.  The
1-forms 
on $\tilde U$
\begin{equation}\label{M88}
\alpha_a = \sum_p F^\alpha_{pa} \theta^p, \quad \beta = 3\sum_p
F^\mu_{\alpha p} \theta^p
\end{equation}
are invariant under such changes of frame, and the bilinear form
\begin{equation}\label{M89}
\psi_a = \alpha_a \otimes \beta
\end{equation}
has the property that $\psi_a(Y_p, Y_p) = 3F^\alpha_{pa} F^\mu_{\alpha
p} = D_{\mu a}$, for every $p$.  By the argument
containing~\eqref{M10}-\eqref{M13}, this implies that
\begin{equation}\label{M91}
D_{\mu a} = 0
\end{equation}
on $\tilde U$, for every $a$.  In the same way, using the second
equation 
in~\eqref{M86}, we conclude that 
\begin{equation}\label{M92}
E_{p\mu}= 0
\end{equation}
on $\tilde U$, for all $p$.  Then~\eqref{M86} implies that
\begin{equation}\label{M93}
F^\mu_{\alpha p} = 0 = F^\mu_{\alpha a}
\end{equation}
on $\tilde U$, for all $a,p$, since $F^\alpha_{pa}$ is never zero on
$\tilde U$, for
all $a,p$.  We remark here for the proof of Corollary~\ref{MC1} below,
that~\eqref{M78} and~\eqref{M93} imply
\begin{equation}\label{M94}
g=h=k=0
\end{equation}
on $\tilde U$, and therefore $\lambda:\tilde U \to \Lambda$ is
reducible, by 
Theorem~\ref{TS2}. For the proof at hand, we use~\eqref{M78} and
~\eqref{M93} in GD(3.39) to get
\begin{equation}\label{M95}
\aligned
0 &= F^\mu_{\alpha aa} = -\frac12 D_{\mu\alpha} \\
0 &= F^\mu_{\alpha a \mu} = D_{a\alpha}
\endaligned
\end{equation}
on $\tilde U$, for all $a,p$, and in GD(3.40) to get
\begin{equation}\label{M96}
0 = F^\mu_{\alpha p \mu} = -D_{p\alpha}
\end{equation}
on $\tilde U$, for all $p$.  Then~\eqref{M85},~\eqref{M91},
\eqref{M92}, \eqref{M95}, and~\eqref{M96} imply that $d\omega^0_0=0$
on $\tilde U$.

In the same way we prove that if $d_1 >0$ on a connected component
$\tilde U$ of $U_2$, then $d\omega^0_0 = 0$ on $\tilde U$.  Therefore,
$d\omega^0_0 = 0$ on all of $U_2$.  If $U_1$ has no interior, then the
closure of $U_2$ in $U$ is all of $U$, and thus $d\omega^0_0=0$ on all
of $U$, by continuity.  If $U_1$ has nonempty interior, then we have
proven that $d\omega^0_0=0$ on this interior, and therefore
$d\omega^0_0=0$ on all of $U$, which is the closure in $U$  of the
union of the
interior of $U_1$ with $U_2$.
\end{proof}

\begin{corollary}\label{MC1} Suppose $m_1=m_2 \geq 2$, $m_3=m_4=1$,
and $r=-1$.  If the Dupin hypersurface is irreducible, then
for any second order frame field $Y:U\to G$ along it, we have
\begin{equation}\label{M97}
d_1+e_2=0
\end{equation}
on $U$.
\end{corollary}

\begin{proof}  Let $Y:U\to G$ be any second order frame field along
the Dupin hypersurface.  We know that $D_1$ and $E_2$ are scalar
matrices at every point of $U$.  Any change of second order frame
field is given by~\eqref{GD(3.14)}, and thus GD(3.32) shows that under
such a change the function $d_1+e_2$ is multiplied by an everywhere
positive function on $U$.  It follows that if~\eqref{M97} holds for
some second order frame field in $U$, then it holds for any second
order frame field in $U$.  

By Proposition~\ref{PM1}, about any point
$x\in U$ there exists a second order frame field for which $d_1 = e_2$.
Seeking a contradiction, suppose that
$d_1(x)\neq 0$, for some $x$ in the domain of this frame field.
Shrinking the domain, if necessary,
we may assume that $d_1 \neq 0$ on the whole 
domain of the frame field, and then the proof of Theorem~\ref{PM3}, as
remarked after~\eqref{M94},
shows that the Dupin hypersurface is reducible on some open subset of
$x$, and thus it is reducible by Proposition~\ref{5P2}.  This
contradicts our assumption that the Dupin hypersurface is 
irreducible.  Hence, $d_1$ must be zero at every point of its
domain.  
\end{proof}

\begin{corollary}\label{MC2} Suppose $m_1=m_2 \geq 1$, $m_3=m_4=1$,
and $r=-1$.  If $Y:U \to G$ is a second order frame field along the
Dupin hypersurface such that one of $\1'$, $\2'$, $\3'$,
or $\4'$ 
is nonempty (see Definition~\ref{DS1}), then the hypersurface is
reducible.
\end{corollary}

\begin{proof} 
If $\1'$ or $\2'$ is nonempty, then there is no loss in generality in
assuming that $\1'$ is nonempty, in which case the result follows from
Theorem~\ref{PU1} and Theorem~\ref{PM3}.  If $\3'$ or $\4'$ is
nonempty, there is no loss in generality in assuming that $\3'$ is
nonempty.  Since $m_3=1$, this means that the functions $f$, $h$, and
$k$, defined in~\eqref{3}, are identically zero on $U$.  Therefore, the
result follows from Theorem~\ref{TS2}.
\end{proof}

\begin{theorem}\label{MT3} Suppose the connected proper Dupin
hypersurface $\lambda: M^{n-1} \to \Lambda^{2n-1}$ has four distinct
curvature spheres with multiplicities
$m_1=m_2 \geq 1$, $m_3=m_4=1$, and Lie curvature $r = -1$.  If
$\lambda$ is irreducible, then it is Lie equivalent to an
isoparametric hypersurface.
\end{theorem}

\begin{proof} 
The case $m_1=m_2=1$ was handled in \cite[Theorem 4.1, p.\ 33]{CJ}.

Assume now that $m_1=m_2 \geq 2$.
In order to apply Theorem~\ref{Th11}, we must find a second order
frame field, defined on a neighborhood of any given point, that
satisfies the conditions of Theorem~\ref{Th11}.  Let $x \in M$ and let
$Y:U \to G$ be a second order frame field, where $x \in U$.  Then
$d_1 +e_2 =0$, by Corollary~\ref{MC1}, and $d\omega^0_0=0$ on $U$, by
Proposition~\ref{PM2}.  By Proposition~\ref{PM1},
we may adjust $Y$ so that
$d_1 = 0 = e_2$ on $U$.  As observed in~\eqref{M16} of the proof of
Proposition~\ref{PM2}, we also have $e_3 = -e_4$ and $d_3 = d_4$.

Consider the functions $A_a$, for any $a \in \1$, defined
in~\eqref{U11a} in the proof of Theorem~\ref{PU1}. 
Irreducibility implies that for each $a \in \1$, $A_a$ must be
positive on a dense open 
subset of $U$, by Theorem~\ref{PU1}.  In the same way, for each $p$,
the analytic function
\begin{equation}\label{M99}
A_p = |v_{p\alpha}|^2 + |v_{p\mu}|^2 = 2\sum_a ((F^\alpha_{pa})^2 +
(F^\mu_{pa})^2) + 8(F^\mu_{\alpha p})^2
\end{equation}
must be positive on a dense open subset of $U$.
There is also a dense open subset of $U$ on which
$E_1$ and $D_2$ can be diagonalized by our frame field.
The intersection of these three dense open subsets is a dense open
subset of $U$.  Let $W$ be a connected component of it.
Then~\eqref{U9} holds on $W$, and $A_a$ is positive on $W$ for all
$a \in \1$, so we have
\begin{equation}\label{M98}
e_a = -e_3
\end{equation}
on $W$, for all $a \in \1$.  Similarly, $A_p$ is positive on $W$ for
all $p$, so by~\eqref{U10} we have
\begin{equation}\label{M100}
d_p = d_3
\end{equation}
on $W$, for all $p \in \2$.  Therefore, on each connected component of
this open dense subset we have $E_1 = -E_3$ and $D_2 = D_3$ are scalar
matrices for some choice of second order frame field, hence for all
choices of second order frame fields.  By continuity, it follows that
this is true on all of $U$ for $Y$, that is,
\begin{equation}\label{M101}
e_1 = -e_3, \quad d_2 = d_3
\end{equation}
on $U$.  Now GD(3.42iii) implies that
\begin{equation}\label{M102}
(m_1 + \frac 12)(d_3 + e_3) = -(m_1 + \frac12) (d_3 + e_3)
\end{equation}
on $U$, which implies that
\begin{equation}\label{M103}
d_3 = -e_3
\end{equation}
on $U$.  Plugging our known values of $d_1,\dots, e_4$
into~\eqref{3.36scalar}, we find that~\eqref{A2i}
holds on $U$.  Finally, it follows now from~\eqref{3.36A1} that
\begin{equation}\label{M104}
|v_{a\alpha}|^2 = \frac12 A_a
\end{equation}
which is positive on a dense open subset of $U$.  Therefore, all the
conditions of Theorem~\ref{Th11} are satisfied and we may conclude
that $\lambda$ is Lie equivalent to an isoparametric hypersurface.
\end{proof}

\begin{remark}\label{Rem34}
There exist reducible proper Dupin hypersurfaces with four curvature
spheres satisfying the conditions $m_1=m_2$, $m_3=m_4$, and $r=-1$
that are not Lie equivalent to an isoparametric hypersurface (see
\cite[pp.\ 107-108]{Cec6}). 
\end{remark}

\section{Compact proper Dupin hypersurfaces}\label{Compact}
Our main goal in this section is to prove the following result.

\begin{theorem}\label{5T1}  Let $W^{d-1}$ be a compact, connected,
proper Dupin hypersurface immersed in $S^d$ (or $\R^d$) with $g >2$
distinct principal curvatures.  Then $W$ is irreducible; that
is, the Legendre submanifold induced by the hypersurface $W$ is
irreducible.
\end{theorem}

\begin{remark}\label{5R1} Thorbergsson \cite{Th1} proved that a
compact proper Dupin hypersurface $W^{d-1}$ immersed in $S^d$ must be
taut, and a taut immersion must be an embedding (see Cecil-Ryan
\cite[p.\ 121]{CecilRyan}), so $W$ must be
embedded in $S^d$.  In this same paper he also proved that the number
$g$ of distinct principal curvatures of a 
compact proper Dupin hypersurface must be 1, 2, 3, 4, or 6,
the same as for an isoparametric hypersurface in
a sphere.
\end{remark}

\begin{proof}  The proof is by contradiction.  We assume that $W^{d-1}
\subset S^d$ 
is reducible and obtain a contradiction.  The proof is essentially the
same as the proof of Theorem~2 of \cite{Cec5}
(see also Theorem~2.11 of \cite[p.\ 148]{Cec6}).  In that
theorem, however, we assumed that the Dupin hypersurface $M^{n-1}
\subset S^n$ to which $W$ is reducible is immersed in $S^n$.  In the
present proof we do not need to make such an assumption because we can
handle the situation as follows.

Assume that there exists a reducible, compact, connected, proper Dupin
hypersurface $W^{d-1}$ immersed in $S^d$ with $g>2$ principal
curvatures.  Let $\nu: W \to \Lambda^{2d-1}$ be the Legendre
submanifold induced by this immersion.  By Proposition~\ref{5P1},
$\nu$ is Lie equivalent to a proper Dupin submanifold $\mu: W \to
\Lambda^{2d-1}$ that is obtained from a lower dimensional proper Dupin
submanifold $\lambda: M^{n-1} \to \Lambda^{2n-1}$ by one of the three
standard constructions of Pinkall \cite{Pinkall}.  As shown in Section~4.2 of \cite{Cec6}, all
three of the constructions begin with a proper Dupin submanifold
$\lambda: M^{n-1} \to \Lambda^{2n-1}$ and produce a Dupin submanifold
$\mu: M \times S^m \to \Lambda^{2(n+m)-1}$.  Thus, $W$ is
diffeomorphic to $M\times S^m$, and $M$ must be compact since $W$ is
compact.  We are also 
assuming that $\nu$ (and thus $\mu$) has $g>2$ distinct curvature
spheres at each point.  As shown in Propositions~2.1, 2.3 and~2.5 of
Section~4.5 of \cite{Cec6}, however, for $\mu$ obtained by the tube
and cylinder constructions, there always exist points on $M\times S^m$
at which the number of distinct curvature spheres is two.  Thus, $\mu$
cannot be obtained by the tube or cylinder construction.

Thus, if $g >2$, then $\mu$ must be obtained from $\lambda$ by the
surface of revolution construction.  Proposition~2.7 of \cite[p.\
144]{Cec6} shows that for the surface of revolution construction, the
number $k$ of distinct curvature spheres on $M$ must be $g-1$ or $g$.
We also have the following relationship between the sum $\beta$ of the
$\Z_2$-Betti numbers of $W$ and $M$,
\begin{equation*}
\beta(W) = \beta(M\times S^m) = 2\beta(M)
\end{equation*}
On the other hand, Thorbergsson \cite{Th1} showed that for a compact
proper Dupin hypersurface embedded in $S^d$, $\beta$ is equal to twice
the number of distinct curvature spheres.  Thus, we have $\beta(W) =
2g$. 

If the point sphere map of $\lambda: M \to \Lambda^{2n-1}$ is an
immersion, or even if there is a Lie sphere transformation $A$ such
that the point sphere map of $A\lambda$ is an immersion, then we have
an immersed proper Dupin hypersurface $f:M \to S^n$ to which
Thorbergsson's theorem applies, and $\beta(M) = 2k$, where $k$ is the
number of distinct curvature spheres of $M$.  Thus, we have
\begin{equation*}
\beta(W) = 2g, \quad \beta(M) = 2k
\end{equation*}
where $k = g-1$ or $k=g$, and
\[
\beta(W) = 2 \beta(M)
\]
Combining these equations, we get
\[
2g = 2(2k) = 4k
\]
for $k = g-1$ or $k=g$.  Clearly, $k=g$ is impossible, and $k=g-1$
yields
\[
2g = 4(g-1) = 4g-4
\]
and thus $g=2$, contradicting the assumption that $g>2$.

It remains to show that in the case of the surface of revolution
construction, 
there is a Lie sphere transformation $A$ such
that the point sphere map of $A\lambda$ is an immersion.  That follows
from the following lemma.

\begin{lemma}\label{5L1} Let $\mu: M^{n-1} \times S^m \to
\Lambda^{2(n+m)-1}$ be a Legendre submanifold which is obtained from a
proper Dupin submanifold $\lambda: M \to \Lambda^{2n-1}$ by the
surface of revolution construction.  If there exists a Lie sphere
transformation $B$ such that the point sphere map of $B\mu$ is an
immersion, then there exists a Lie sphere transformation $A$ such that
the point sphere map of $A\lambda$ is an immersion.
\end{lemma}

We shall first apply this lemma to complete the proof of
Theorem~\ref{5T1} and then give the proof of the lemma.  Since the
point sphere map of $\nu: W \to \Lambda^{2d-1}$ is given to be an
immersion , and $\nu$ is Lie equivalent to $\mu: W \to
\Lambda^{2d-1}$, we know that the Lie sphere transformation $B$ in the
lemma exists.  Therefore, a Lie sphere transformation $A$ exists such
that the point sphere map of $A\lambda$ is an immersion, which is what
we need to complete the proof of Theorem~\ref{5T1}.
\end{proof}

\begin{proof}[Proof of Lemma~\ref{5L1}] We begin by reviewing the
surface of 
revolution construction (see \cite[pp.\ 141-144]{Cec6}).  Let
$e_0,\dots,e_{n+m+2}$ be an orthonormal basis of $\R^{n+m+3}_2$ with
$e_0$ and $e_{n+m+2}$ timelike and the rest spacelike.  Let
$\rp^{n+m+2}$ be the projective space determined by $\R^{n+m+3}_2$
with corresponding Lie quadric $Q^{n+m+1}$.  Let
\begin{equation}\label{5e2}
\R^{n+3}_2 = \span\{e_0,e_1,\dots,e_{n+1},e_{n+m+2}\} \subset
\R^{n+m+3}_2
\end{equation}
and let $\rp^{n+2}$ and $Q^{n+1}$ be the corresponding projective
space and Lie quadric.  Let $\Lambda^{2n-1}$ and $\Lambda^{2(n+m)-1}$
be the space of projective lines on $Q^{n+1}$ and $Q^{n+m+1}$,
respectively.  Finally, let $u_k = e_{k+1}$ for $k = 1,\dots,n+m$, and
let
\begin{equation}\label{5e3}
\R^n = \span \{ e_2,\dots,e_{n+1} \} = \span \{u_1,\dots,u_n \}
\end{equation}
\begin{equation}\label{5e4}
\R^{n+m} = \span \{e_2,\dots,e_{n+m+1}\} = \span \{ u_1,\dots,u_{n+m}
\}
\end{equation}
Consider a proper Dupin submanifold $\lambda: M^{n-1} \to
\Lambda^{2n-1}$ with $g$ distinct curvature spheres.  We can
parametrize $\lambda$ by using the Euclidean projection $f:M \to \R^n$
and Euclidean field of unit normals $\xi: M \to \R^n$ as follows (see
\cite[p.\ 82]{Cec6}),
\begin{equation}\label{5e5}
\lambda = [k_1,k_2]
\end{equation}
where
\begin{equation}\label{5e6}
k_1 = (1+f\cdot f, 1-f\cdot f, 2f, 0)/2,\quad k_2 = (f\cdot \xi,
-f\cdot \xi, \xi, 1)
\end{equation}
The map $[k_1]:M \to Q^{n+1}$ is the point sphere map of $\lambda$,
and the map $[k_2]: M \to Q^{n+1}$ is the tangent hyperplane map of
$\lambda$. 

We want to construct a Legendre submanifold $\mu$ by ``revolving'' $f$
around an ``axis'' $\R^{n-1} \subset \R^n \subset \R^{n+m}$, for
$\R^{n+m}$ as in~\eqref{5e4}.  The domain of $\mu$ will be $M \times
S^m$.  Note that we do not assume that $f$ is an immersion, nor that
the range of $f$ is disjoint from the axis $\R^{n-1}$.

For simplicity, we assume that the axis $\R^{n-1}$ contains the origin
of $\R^n$ and that the orthonormal basis vectors have been chosen so
that
\begin{equation}\label{5e7}
\R^{n-1} = \span \{u_1,\dots, u_{n-1}\} \subset \R^n \subset \R^{n+m}
\end{equation}
for $\R^n$ and $\R^{n+m}$ as in~\eqref{5e3} and~\eqref{5e4},
respectively.  We write the sphere $S^m$ in the form
\begin{equation}\label{5e8}
S^m = \{ y = y_0 u_n + \dots + y_m u_{n+m} : y_0^2 + \dots + y_m^2 = 1
\}
\end{equation}
We now define a new Legendre submanifold $\mu:M \times S^m \to
\Lambda^{2(n+m)-1}$ by its Euclidean projection $F:M \times S^m \to
\R^{n+m}$ and its Euclidean field of unit normals $\eta: M\times S^m
\to \R^{n+m}$, so that $\mu = [K_1,K_2]$, where
\begin{equation}\label{5e9}
K_1 = (1+F\cdot F, 1-F\cdot F, 2F, 0)/2, \quad K_2 = (F\cdot \eta,
-F\cdot \eta, \eta,1)
\end{equation}
and the maps $F$ and $\eta$ are defined as follows.  First we
decompose the maps $f$ and $\xi$ into components along $\R^{n-1}$ and
orthogonal to $\R^{n-1}$ in $\R^n$ and write,
\begin{equation}\label{5e10}
f(x) = \hat f(x) + f_n(x) u_n, \quad \hat f(x) \in \R^{n-1}
\end{equation}
\begin{equation}\label{5e11}
\xi(x) = \hat\xi(x) + \xi_n(x) u_n, \quad \hat\xi(x) \in \R^{n-1}
\end{equation}
Then for $x \in M$, $y \in S^m$, we define the maps $F$ and $\eta$
in~\eqref{5e9} by
\begin{equation}\label{5e12}
F(x,y) = \hat f(x) + f_n(x) y
\end{equation}
\begin{equation}\label{5e13}
\eta(x,y) = \hat \xi(x) + \xi_n(x) y
\end{equation}
In \cite[pp.\ 141-144]{Cec6}, we show that the curvature spheres of
$\mu$ are of two types.  The first type is
\begin{equation}\label{5e14}
[K(x,y)] = [\xi_n(x)K_1(x,y) - f_n(x) K_2(x,y)]
\end{equation}
This is the new curvature sphere introduced by the surface of
revolution construction.  The second type is
\begin{equation}\label{5e15}
[K(x,y)] = [rK_1(x,y) + sK_2(x,y)]
\end{equation}
where $r,s$ are real numbers such that
\begin{equation}\label{5e16}
[k(x)] = [rk_1(x) + sk_2(x)]
\end{equation}
is a curvature sphere of $\lambda$ at $x\in M$.

Now we turn to the question of whether there exists a Lie sphere
transformation $A$ such that the point sphere map of the Legendre
submanifold $A\lambda: M \to \Lambda^{2n-1}$ is an immersion.  The
point sphere map $Z(x)$ of the Legendre submanifold $A\lambda$ is
determined by the condition
\begin{equation}\label{5e17}
\langle Z(x), e_{n+m+2}\rangle =0
\end{equation}
The map $Z(x)$ is an immersion if and only if $Z(x)$ is not a
curvature sphere of $A\lambda$, for any $x\in M$.  Since the curvature
spheres of $A\lambda$ are of the form $A[k]$, where $[k]$ is a
curvature sphere of $\lambda$, we see that $Z$ is an immersion if and
only if
\begin{equation}\label{5e18}
\langle Ak(x), e_{n+m+2} \rangle \neq 0
\end{equation}
for all curvature spheres $[k(x)]$ of $\lambda$ at all points $x\in
M$.  If we apply the Lie sphere transformation $A^{-1} \in O(n+1,2)$
to~\eqref{5e18}, we get the condition
\begin{equation}\label{5e19}
\langle k(x), A^{-1}e_{n+m+2} \rangle \neq 0
\end{equation}
for all curvature spheres $[k(x)]$ of $\lambda$ at all points $x \in
M$.  Note that $v = A^{-1} e_{n+m+2}$ is a unit timelike vector in
$\R^{n+3}_2$.  Thus, there exists a Lie sphere transformation $A$ such
that the point sphere map $Z$ of $A\lambda$ is an immersion if and
only if there exists a unit timelike vector $v \in \R^{n+3}_2$ such
that
\begin{equation}\label{5e20}
\langle k,v\rangle \neq 0
\end{equation}
for all curvature sphere maps $k$ of $\lambda$.

It is a hypothesis of Lemma~\ref{5L1} that there exists a Lie sphere
transformation $B \in O(n+m+1,2)$ such that the point sphere map $W$
of the Legendre submanifold $B\mu$ is an immersion.  Thus, there
exists a unit timelike vector $q \in \R^{n+m+3}_2$ such that
\begin{equation}\label{5e21}
\langle K,q \rangle \neq 0
\end{equation}
for all curvature sphere maps $K$ of $\mu$.  We can write $q$ in
coordinates as
\begin{equation}\label{5e22}
q = (q_0,q_1,\hat q,w,q_{n+m+2})
\end{equation}
where $\hat q = (q_2,\dots,q_n)\in \R^{n-1}$ and $w =
(q_{n+1},\dots,q_{n+m+1})$.  For a curvature sphere $K(x,y)$ as
in~\eqref{5e15}, one can compute that
\begin{equation}\label{5e23}
\aligned
\langle &K(x,y), q \rangle = -q_0(r(\frac{1+f\cdot f}2) + sf\cdot \xi )
\\
&+ q_1(r(\frac{1-f\cdot f}2) - sf\cdot \xi) + (r\hat f(x) + s\hat
\xi(x)) \cdot \hat q \\
&+ (r f_n(x) + s\xi_n(x))(y\cdot w) - s q_{n+m+2}
\endaligned
\end{equation}
since $F\cdot F = f\cdot f$, $\eta \cdot \eta = \xi \cdot \xi$, and
$F\cdot \eta = f \cdot \xi$.
Note that all terms depend only on $x \in M$ except the term
\begin{equation}\label{5e24}
(rf_n(x) + s\xi_n(x))(y\cdot w)
\end{equation}
If we take $y = u_n$, then
\begin{equation}\label{5e25}
\langle K(x,y),q \rangle = \langle k(x),v \rangle
\end{equation}
for $[k(x)]$ as in~\eqref{5e16} and $v = \pi(q)$, where $\pi$ is
orthogonal projection of $\R^{n+m+2}_2$ onto its subspace $\R^{n+3}_2$
in~\eqref{5e2} given by
\begin{equation}\label{5e26}
v = \pi(\sum_{i=0}^{n+m+2} q_i e_i) = \sum_{i=0}^{n+1} q_i e_i +
q_{n+m+2} e_{n+m+2}
\end{equation}
Note that $v$ is timelike, since $q$ is a unit timelike vector, and
\begin{equation}\label{5e27}
\langle v,v \rangle = \langle q, q \rangle - (q_{n+2}^2 + \dots +
q_{n+m+1}^2)
\end{equation}
Thus, $v$ is a timelike vector such that 
\begin{equation}\label{5e28}
\langle k,v \rangle \neq 0
\end{equation}
for all curvature sphere maps $[k]$ of $\lambda$.  So there exists a
Lie sphere transformation $A \in O(n+1,2)$ such that the point sphere
map of $A\lambda$ is an immersion.
\end{proof}

\begin{theorem}\label{5T2} Let $M$ be a compact connected proper Dupin
hypersurface immersed in $R^n$ with four distinct principal curvatures
having multiplicities $m_1=m_2 \geq 1$, $m_3 = m_4 = 1$, and 
constant Lie curvature.  Then $M$
is Lie equivalent to an isoparametric hypersurface.
\end{theorem}

\begin{proof} As noted in Remark~\ref{5R1}, $M$ must, in fact, be
embedded in $\R^n$. 
Miyaoka \cite[p.\ 252]{Mi2} showed that if $M$ is a
compact connected proper Dupin hypersurface embedded in $\R^n$ with
four distinct principal curvatures and constant Lie curvature $r$,
then $r$ must equal $-1$.  (Miyaoka's therorem states that $r = 1/2$,
but in that case the order of the principal curvatures can be
rearranged so that $r = -1$).  Then
Theorem~\ref{5T1} implies that $M$ is irreducible, and then
Theorem~\ref{MT3} implies that $M$ is Lie equivalent to an
isoparametric hypersurface.
\end{proof}

\begin{remark}\label{5R2} Theorem~\ref{5T1} and Corollary~\ref{5C1}
imply that a compact proper Dupin hypersurface immersed in $\R^n$ with
$g=3$ principal curvatures must be Lie equivalent to an isoparametric
hypersurface.  This was first proven by Miyaoka \cite{Mi1}, who used
different methods and did not focus on the notion of irreducibility.
\end{remark}

\bibliography{Bibliography}

\begin{thebibliography}{10}
\providecommand{\url}[1]{{#1}}
\providecommand{\urlprefix}{URL }
\expandafter\ifx\csname urlstyle\endcsname\relax
  \providecommand{\doi}[1]{DOI~\discretionary{}{}{}#1}\else
  \providecommand{\doi}{DOI~\discretionary{}{}{}\begingroup
  \urlstyle{rm}\Url}\fi

\bibitem{Car2}
Cartan, E.: {Sur des familles remarquables d'hypersurfaces
  isopa\-ra\-m\'etriques dans les espaces sph\'eriques}.
\newblock Math. Z. \textbf{45}, 335 -- 367 (1939)

\bibitem{Cec5}
Cecil, T.E.: Reducible {D}upin submanifolds.
\newblock Geom. Dedicata \textbf{32}(3), 281--300 (1989)

\bibitem{Cec6}
Cecil, T.E.: Lie sphere geometry.
\newblock Universitext. Springer-Verlag, New York (1992)

\bibitem{AlgDupin}
Cecil, T.E., Chi, Q.S., Jensen, G.R.: Algebraic properties of {D}upin
  hypersurfaces (2005).
\newblock Preprint

\bibitem{CJ}
Cecil, T.E., Jensen, G.R.: Dupin hypersurfaces with three principal curvatures.
\newblock Invent. Math. \textbf{132}(1), 121--178 (1998)

\bibitem{CJ2}
Cecil, T.E., Jensen, G.R.: Dupin hypersurfaces with four principal curvatures.
\newblock Geom. Dedicata \textbf{79}(1), 1--49 (2000)

\bibitem{CecilRyan2}
Cecil, T.E., Ryan, P.J.: Focal sets, taut embeddings and the cyclides of
  {D}upin.
\newblock Math. Ann. \textbf{236}(2), 177--190 (1978)

\bibitem{CecilRyan}
Cecil, T.E., Ryan, P.J.: Tight and taut immersions of manifolds, \emph{Research
  Notes in Mathematics}, vol. 107.
\newblock Pitman (Advanced Publishing Program), Boston, MA (1985)

\bibitem{GH}
Grove, K., Halperin, S.: Dupin hypersurfaces, group actions and the double
  mapping cylinder.
\newblock J. Differential Geom. \textbf{26}(3), 429--459 (1987)

\bibitem{Mi1}
Miyaoka, R.: Compact {D}upin hypersurfaces with three principal curvatures.
\newblock Math. Z. \textbf{187}(4), 433--452 (1984)

\bibitem{Mi2}
Miyaoka, R.: Dupin hypersurfaces and a {L}ie invariant.
\newblock Kodai Math. J. \textbf{12}(2), 228--256 (1989)

\bibitem{Mi3}
Miyaoka, R.: Dupin hypersurfaces with six principal curvatures.
\newblock Kodai Math. J. \textbf{12}(3), 308--315 (1989)

\bibitem{MiOz}
Miyaoka, R., Ozawa, T.: Construction of taut embeddings and {C}ecil-{R}yan
  conjecture.
\newblock In: Geometry of manifolds (Matsumoto, 1988), \emph{Perspect. Math.},
  vol.~8, pp. 181--189. Academic Press, Boston, MA (1989)

\bibitem{Mun1}
M{\"u}nzner, H.F.: Isoparametrische {H}yperfl\"achen in {S}ph\"aren.
\newblock Math. Ann. \textbf{251}(1), 57--71 (1980)

\bibitem{Mun2}
M{\"u}nzner, H.F.: Isoparametrische {H}yperfl\"achen in {S}ph\"aren. {II}.
  \"{U}ber die {Z}erlegung der {S}ph\"are in {B}allb\"undel.
\newblock Math. Ann. \textbf{256}(2), 215--232 (1981)

\bibitem{Pinkall}
Pinkall, U.: {D}upin hypersurfaces.
\newblock Math. Ann. \textbf{270}, 427 -- 440 (1985)

\bibitem{Pink}
Pinkall, U.: Dupinsche {H}yperfl\"achen in {$E\sp 4$}.
\newblock Manuscripta Math. \textbf{51}(1-3), 89--119 (1985)

\bibitem{PT}
Pinkall, U., Thorbergsson, G.: Deformations of {D}upin hypersurfaces.
\newblock Proc. Amer. Math. Soc. \textbf{107}(4), 1037--1043 (1989)

\bibitem{Singley}
Singley, D.H.: Smoothness theorems for the principal curvatures and principal
  vectors of a hypersurface.
\newblock Rocky Mountain J. Math. \textbf{5}, 135--144 (1975)

\bibitem{Stolz}
Stolz, S.: Multiplicities of {D}upin hypersurfaces.
\newblock Invent. Math. \textbf{138}(2), 253--279 (1999)

\bibitem{Th1}
Thorbergsson, G.: Dupin hypersurfaces.
\newblock Bull. London Math. Soc. \textbf{15}(5), 493--498 (1983)

\end{thebibliography}
\bibliographystyle{spmpsci}

\enddocument
\end